\documentclass[preprint, authoryear]{elsarticle}

\usepackage[T1]{fontenc}
\usepackage[table,xcdraw]{xcolor}
\usepackage{lmodern}
\usepackage{amssymb}
\usepackage{amsmath} 
\usepackage{dsfont}
\usepackage{multirow}
\usepackage{subcaption}
\usepackage{bm}
\usepackage{natbib}
\usepackage{apalike}
\setcitestyle{round}
\usepackage[linesnumbered,ruled,vlined]{algorithm2e}
\usepackage{float}
\usepackage{mathtools}
\usepackage {xcolor}
\usepackage{graphicx}
\usepackage{subcaption}
\usepackage{geometry}
\usepackage{microtype}
\usepackage{mathrsfs}
\usepackage{hyperref}
\hypersetup{
    colorlinks=true,
    linkcolor=blue,
    filecolor=magenta, 
    citecolor=blue,     
    urlcolor=cyan,
    pdftitle={Overleaf Example},
    pdfpagemode=FullScreen,
    }

\begin{document} 

\begin{frontmatter}

\title{Data-driven Optimization for Drone Delivery Service Planning\\with Online Demand}
\author[1]{Aditya Paul}
\ead{aditya.paul@unsw.edu.au}
\author[3]{Michael W. Levin}
\ead{mlevin@umn.edu}
\author[2]{S. Travis Waller}
\ead{s.travis.waller@gmail.com}
\author[4]{David Rey\corref{cor1}}
\ead{david.rey@skema.edu}

\cortext[cor1]{Corresponding author}

\affiliation[1]{organization={School of Civil and Environmental Engineering},
addressline={UNSW Sydney},
postcode={NSW, 2052},
city={Sydney},
country={Australia}}

\affiliation[3]{organization={Department of Civil, Environmental, and Geo-Engineering, University of Minnesota},
postcode={55455},
city={Minneapolis},
country={USA}}

\affiliation[2]{organization={Lighthouse Professorship ‘‘Transport Modelling and Simulation’’, Faculty of Transport and Traffic Sciences, Technische Universität},
city={Dresden},
country={Germany}}

\affiliation[4]{organization={SKEMA Business School, Université Côte d’Azur},
addressline={Sophia Antipolis},
country={France}}

\begin{abstract}
In this study, we develop an innovative data-driven optimization approach to solve the drone delivery service planning problem with online demand. Drone-based logistics are expected to improve operations by enhancing flexibility and reducing congestion effects induced by last-mile deliveries. With rising digitalization and urbanization, however, logistics service providers are constantly grappling with the challenge of uncertain real-time demand. This study investigates the problem of planning drone delivery service through an urban air traffic network to fulfil online and stochastic demand. Customer requests – if accepted – generate profit and are serviced by individual drone flights as per request origins, destinations and time windows. We cast this stochastic optimization problem as a Markov decision process. We present a novel data-driven optimization approach which generates predictive prescriptions of parameters of a surrogate optimization formulation. Our solution method consists of synthesizing training data via lookahead simulations to train a supervised machine learning model for predicting relative link priority based on the state of the network. This knowledge is then leveraged to selectively create weighted reserve capacity in the network and via a surrogate objective function that controls the trade-off between reserve capacity and profit maximization to maximize the cumulative profit earned. Using numerical experiments based on benchmarking transportation networks, the resulting data-driven optimization policy is shown to outperform a myopic policy. Sensitivity analyses on learning parameters reveal insights into the design of efficient policies for drone delivery service planning with online demand.
\end{abstract}


\end{frontmatter}

\section{Introduction}

Unmanned Aerial Vehicles (UAVs) have been growing in both prevalence and application. UAVs provide the advantage of performing a task aerially with limited or no human supervision. Additionally, most UAV designs are light-weight and can be deployed in large numbers. These features cause UAVs to be in high demand in many different domains, particularly logistics. UAV applications broadly fall into two categories: either some form of data acquisition and transmission (also called monitoring or reconnaissance), or payload pickup and delivery. This study falls into the second category and addresses the drone delivery service planning problem (DDSP) under stochastic online demand. The DDSP takes the perspective of a transportation network manager that plans service for drone deliveries. We consider discrete time steps grouped into time intervals. At each time step, customers may submit delivery requests that are to be served by a fleet of UAVs subject to time windows and link capacity constraints. The probability distributions behind request parameters, and the request arrival process, are unknown to the network manager and estimated from historical data. Decisions are made at every time interval to accept or reject the requests that arrive and to route accepted requests. We explore the potential of data-driven optimization methods for this online DDSP problem.
\\
\\
In stochastic optimization problems, some or all parameters are subject to uncertainty. These parameters are not deterministic, but instead emerge from probability distributions or uncertainty sets. Traditional solution methods to stochastic problems exploit \textit{a priori} knowledge of these probability distributions. Stochastic problems can further be sub-divided into two categories: (a) static, where all problem parameters are available beforehand, and some problem specific condition determines if and when the stochastic parameters are instantiated; (b) dynamic, where problem parameters data arrives is revealed over time. Besides, online optimization makes no assumption on the nature of uncertain data and focuses on the immediacy and real-time nature of incoming problem parameters. Online problems are always dynamic in nature. Nevertheless, in many cases, it is reasonable to presume that the probability distribution, or an estimate thereof, can be extracted either through historical data or predictive models \citep{van2010online}. From a computational standpoint, stochastic and online optimization problems lead to significant challenges. The rise of e-commerce has compelled logistic service providers to tackle online demand. In this study, we consider the setting of the static DDSP in the context of online demand, i.e. user requests that are revealed over time.
\\
\\
Our solution methodology is based on an innovative data-driven optimization and lies at the intersection of predictive analytics (machine learning) and prescriptive analytics (stochastic optimization). The literature on the combination of machine learning (ML) and stochastic optimization reveals two distinct avenues for a successful collaboration:

\begin{itemize}
\item{ML-enabled optimization: Combinatorial optimization problems can be very hard to solve, and special techniques or methods can be used to solve a particular kind of problem or problem instance. Instead of relying on experts to create such techniques, algorithms and models can be developed to automatically learn from implicit distribution of problem instances and adapt to those instances. ML has been found to be useful here.}
\item{Prediction and optimization: Many real-world practical applications involve uncertain parameters which directly or indirectly affect the objective function, but which must first be predicted before optimization. In such predict-then-optimize situations, ML is being recognized as a vital tool to make accurate data-driven predictions.}
\end{itemize}

Research into transportation service planning under online demand, with a dedicated integration of ML and mathematical programming (MP), is still in its nascent stages. In this study, we incorporate strategies from both these avenues to develop a customized data-driven optimization approach for the DDSP with online demand. This study makes the following contributions to the field. 

\begin{enumerate}
\item{We extend the integer linear programming (ILP) formulation of the deterministic DDSP to incorporate online demand, i.e. delivery requests are revealed over time and routing decisions must be taken to accommodate accepted requests, and cast this problem as a Markov Decision Process (MDP).}

\item{We propose a novel solution method for the online DDSP where predictive prescriptions are computed over time. Historical demand data including spatial and temporal distributions, as well as the request arrival process are used to determine link-based features that inform a surrogate optimization formulation.}

\item{We demonstrate through numerical experiments that our data-driven optimization approach outperforms a myopic optimization approach. Sensitivity analysis reveal how learning parameters affect optimal policies for the online DDSP.} 

\end{enumerate}

Since we base our solution method on link-based features, our approach is applicable to other network-based logistics systems. The rest of the paper is organized as follows. In Section~\ref{lit}, we review the relevant literature. In Section~\ref{prob_des}, we present the online DDSP formulation and outline the MDP framework. Section~\ref{sol_met} develops the solution methodology and the proposed algorithms. Section~\ref{num_exp} reports the numerical experiments and their results. Our findings and research perspectives are discussed in Section~\ref{con}.

\section{Literature review}
\label{lit}

We examine the state-of-the art relevant to the online DDSP. Section~\ref{pred_and_presc} reviews recent advances in data-driven optimization, notably the interplay between ML and optimization. Section~\ref{online_intro} records the advent of stochastic and online demand as well as the latest solution approaches arising in urban logistics. Finally, Section~\ref{uav_intro} outlines the widespread adoption of UAVs in commercial and industrial sectors.

\subsection{Data-driven optimization}
\label{pred_and_presc}

Machine learning (ML) has proved useful for optimization purposes along two lines: to assist directly within the optimization process by exploiting particular problem characteristics, or as a data-driven predictive framework to estimate uncertain parameters that are involved in the optimization formulation. We review these two streams of literature hereafter.

\subsubsection{Machine learning-enabled optimization}
\label{ml-opt}

It is common in optimization that a solution algorithm is tuned to a specific problem scenario by refining its parameters or policy to that scenario. Such tuning (or adaptation) need not come from research expertise and intuition. ML can assist an optimization algorithm to this end~\citep{bengio2021machine}. This direct assistance provided by ML can be analysed based on its motivation or on its structure. In terms of motivation, ML can either: (a) replace some heavy computation with a fast, reliable approximation or (b) explore the decision space (defined by the optimization model) and discover better policies. As for structure, the collaboration between ML and optimization can be sequential, whereby the ML model supports the optimization algorithm with valuable information. Or it can be parallel, wherein the ML and MP components work alongside each other. If one of the components works within the framework of the other, then the solutions of the nested component help the external one make better decisions.
\\
\\
\cite{kruber2017learning} provide an illustration of the sequential implementation. They use supervised ML to decide if a given Dantzig-Wolf decomposition will be effective in solving instances of a mixed-integer linear programming (MILP) problem. \cite{bonami2018learning} address a similar question. They employ supervised ML to determine if linearizing a given instance of a mixed integer quadratic programming problem will reduce solve time. An example of parallel implementation is found in \cite{lodi2017learning}. The ML model works within a branch-and-bound framework, and is used to make decisions about branching. Specifically, the goal of the ML component is to learn an efficient policy for selecting the branching variable and branching direction. \cite{ yilmaz2024expandable} address sequential combinatorial optimization problems in industrial contexts, where problem instances arise with slight parameter variations. A neural network generates a relaxed problem instance by predicting binding constraints, then uses this relaxed instance to eliminate infeasible predictions of decision variables through rapid feasibility checks. This reduces solution time by up to three orders of magnitude, while maintaining an optimality gap below 0.1\%. \cite{larsen2022predicting} introduce a supervised ML model to quickly predict expected tactical descriptions of operational solutions in stochastic optimization. This addresses a common challenge in two-stage stochastic programming, where the second stage is computationally demanding. Their approach aids in solving the overall two-stage problem by circumventing the need for online generation of multiple second stage scenarios and solutions. Another implementation possibility is to leverage a MP formulation to aid an ML algorithm in the solution process. \cite{tran2012knapsack} employ an unbounded knapsack formulation to improve the performance of the $\epsilon$-greedy algorithm for multi-armed bandits, where the knapsack solution governs the probability of pulling a given arm of the bandit. A comprehensive overview of ML-enabled optimization methods is reported in~\cite{bengio2021machine}.

\subsubsection{Predict-then-optimize}
\label{pto}

Many decision-making problems deal with both prediction and optimization (or decision). These two tasks are complex on their own and are hence often treated separately and sequentially. The prediction model predicts the unknown parameters of the optimization model, which the latter uses to perform the optimization. ML has been a popular candidate for the prediction task in recent times. But the conventional approach is to train the ML prediction model to minimize the prediction error. This does not consider the nature and attributes of the successive optimization problem and may result in sub-optimal decisions~\citep{elmachtoub2022smart}. As a brief but standard example, consider two possible paths between an origin O and a destination D. Travel times for both paths are unknown to the decision-maker, but may be predicted with contextual data (weather, congestion, etc.). Suppose the actual travel times are 10 and 20 minutes each, respectively. ML model 1 predicts travel times on the two paths to be 15 and 13 minutes (respectively), while ML model 2 predicts the same to be 30 and 40 minutes. Although model 1 has better accuracy, it causes a worse decision.
\\
\\
Two frameworks have been proposed which integrate prediction with optimization by taking the downstream optimization model into account \citep{yan2022integrating}. The first framework is called smart “predict, then optimize” (SPO), and it evaluates the performance of the ML prediction model based on the final decision error induced by the prediction made, not on the initial prediction error \citep{elmachtoub2022smart}. The authors propose a design which applies to the optimization of a linear objective functions over a convex feasible region. They introduce the SPO loss function which captures the decision error of a prediction by utilizing the objective function of the downstream optimization model. Since training with the SPO loss function can be cumbersome, a more manageable surrogate loss function called SPO+ loss is also introduced, which is equivalent to the original loss under specific conditions. The second framework, predictive prescriptions, aims to harness the joint distributions between decision variables and auxiliary data, both of which are uncertain \citep{bertsimas2020predictive}. The principle of predictive prescriptions is to simulate various scenarios in which the values of the unknown parameters are estimated by the ML prediction model. These estimations lead to a cost in the optimization process, and the costs are weighted by the data-driven exploration of the scenarios. The optimal decision policy then minimizes the weighted sum costs thus generated. Both SPO and predictive prescriptions framework have been increasingly adopted by the scientific community and a survey of contextual optimization methods for decision-making under uncertainty is available in~\cite{sadana2023survey}.

\subsection{Stochastic and online demand}
\label{online_intro}

Online demand, sometimes referred to as on-demand logistics, deals with that aspect of e-commerce which provides an affordable delivery service to satisfy the customer's place and time requirements at short notice. It concerns the sale of immediately deliverable products, sometimes at the cost of short-term losses incurred by the e-commerce business or logistic service provider. With the consolidation of the requisite digital mechanisms (internet ubiquity, personal computers, mobile phones, online payment apps) and supply chain infrastructure (warehouses, highways, enterprise resource planning, GPS), e-commerce has experienced a consistent surge in online demand in the 21st century. E-commerce is now experiencing over $10\%$ growth per annum worldwide, and B2C e-commerce deliveries rose by nearly $25\%$ in 2020 alone \citep{lozzi2022demand}. In addition, the disruption caused by COVID-19 has caused an acceleration in internet usage and a lowering of traditional offline purchases, both factors contributing further to e-commerce home deliveries~\citep{han2022covid}.
\\
\\
Even in the 2000s the research literature observed the benefits of transitioning from brick-and-mortar retail to online e-commerce. Particularly, engaging in both physical offline sales and e-commerce -- what is called multi-channel retailing, or click-and-mortar -- was identified as a crucial paradigm shift. In the USA, more than $50\%$ of multi-channel retailers reported positive operating margins as early as 2001 \citep{kim2005consumer}. Investing in e-commerce gave companies access to new markets at relatively low capital costs. Those customers who shopped both online and offline exhibited greater loyalty to the business. The stable customer base and retailer trust acquired from years of brick-and-mortar operation could be extended to the online domain with e-commerce expansions. Engaging with online demand also gives businesses low-cost direct contact with their customers \citep{fruhling2000impact}. This not only helps build customer relationships faster, but also enables the construction of customer profiles, the analysis of customer interests and purchasing behaviour, and comparison of similar profiles to make accurate and relevant recommendations. The comprehensive and coherent advantages of online e-commerce have even led to the evolution of omni-channel retailing, which is the modification of multi-channel retailing such that customers can experience and exploit a seamless array of multiple purchase possibilities \citep{verhoef2015multi}. 
\\
\\
We next discuss the impact of online demand on vehicle routing, both in terms of emerging obstacles and recent solution approaches. We then describe solution strategies which incorporate ML, in line with the methodologies outlined in Section \ref{pred_and_presc}.

\subsubsection{Online demand in vehicle routing}

The benefits of e-commerce come with the challenges of online demand. By its dynamic and uncertain nature, online demand is difficult to fulfill. Failed first deliveries can account for up to $60\%$ of all orders for a service provider \citep{song2009addressing}, the primary reason being the absence of the customer during package handover. This is partially due to modern lifestyle changes such as increase in single person households and flexible (but demanding) work patterns. But a major contributing factor also is that no delivery time slot is explicitly agreed upon between the customer and the e-commerce retailer. This leads to several detrimental effects. If packages are left unattended then lack of security is a concern. Repeated deliveries or customers personally travelling to collection points result in economic losses as well as excess carbon emissions \citep{edwards2010carbon}. A high rate of unsatisfied or partially satisfied demand also strains driver-customer relationships \citep{masorgo2023you}. However, if a delivery time window is agreed upon through some form of online demand management, then the success of last mile deliveries is shown to increase \citep{van2016improving}.
\\
\\
Many studies dealing with online demand employ some form of heuristic approach to modify routes as and when new demand requests emerge. \cite{hong2019routing} consider a shipping company that faces online demand, but delivers products to pre-determined collection centres that the customer must visit on their own. They implement an ant colony based heuristic which changes the previously decided collection centre for a particular customer in light of incoming demand requests. Mobility services also fall under stochastic and online demand, where the commodity sold is the route itself. \cite{bertsimas2019online} analyse online taxi ride-sharing with time windows for pickups. They employ a cost metric to measure waiting or empty driving time between customers. The network is pruned by keeping, for each customer, the top $k$ subsequent customers with the lowest costs. A maxflow heuristic is applied by reducing pickup time windows to a single randomly selected time instant. The resulting route arcs are retained, and the process is iteratively repeated until a specified number of arcs is reached. Finally, a mixed-integer formulation is applied to the modified network for obtaining results. Consolidation of delivery requests can also save costs by revealing efficient routes. \cite{munoz2024large} explore stochastic demand where customer requests arrive daily. The city of operation is divided into regions, each with an expected customer density, estimated average velocity and travel distance per order. A region-specific postponement cost function is proposed for quantifying the delay of a request in a given region to a later date. A tabu-search metaheuristic modifies vehicle routes by identifying better (lower cost) routes in their local neighbourhood. The neighbourhood is traversed by applying insertion, swap and inversion operators that transform the customer sequence of a vehicle route. 

\subsubsection{ML-based solution methods}

Using ML to assist optimization was given much attention in the transport literature. The same-day delivery problem is a rich base for such joint implementations as mentioned in Section \ref{ml-opt}. \cite{chen2023same} study same-day delivery under fairness, dividing the service area into regions. Service availability is defined as the ratio of expected accepted requests per day to the expected overall number of requests. The objective is to maximize both the service availability in the entire area and the minimum service availability across all regions. The learning agent receives rewards for each fulfilled request, guiding it to improve its policy by considering changes in both components of the objective function. \cite{feng2022coordinating} apply reinforcement learning to the real-time ride-sourcing problem with incorporation of public transport services. For every customer ride request, feasible routes are heuristically generated. For every driver capable of responding to the request, the learning agent estimates the perceived value of the route-driver pairs. Finally an integer linear program is solved to determine the optimal allocation of drivers to routes. If electric vehicles are involved, then charge depletion adds an additional factor of uncertainty. \cite{basso2022dynamic} study the electric vehicle routing problem where both arrival of customer requests and energy consumption per unit road length are stochastic. Charging locations are deterministic. At each state of the MDP formulation, the reinforcement learning agent estimates the energy cost of moving to possible subsequent states and, in doing so, the probability of battery charge falling below a safety threshold. 
\\
\\
The predict-then-optimize frameworks discussed in Section~\ref{pto} have also been applied in stochastic transport optimization. \cite{chu2021data} implement the SPO method to solve the last-mile problem for an online food delivery service. They consider driver behaviour, traffic conditions and route lengths as input features to predict travel times. Simulated annealing is used to generate feasible routes, which are further updated with customer exchange operators to derive new routes. \cite{baty2024combinatorial} design a machine learning pipeline which relies on a downstream optimization component to improve predictions. The authors study the dynamic vehicle routing problem with time windows (VRPTW), where all requests must be fulfilled and fleet size is unlimited. They demonstrate that solving an equivalent prize-collecting VRPTW, where each unfulfilled customer is assigned a prize or profit, addresses this version of the VRP. An ML model determines the optimal allocation of prizes, yielding a near-optimal solution for the original problem. A hybrid genetic search algorithm solves the prize-collecting VRPTW based on input from the ML component. The ML model refines itself through a loss function, which measures the disparity between the target and acquired objective function values. The target objective value is derived by solving a static VRPTW from historical data, utilizing the same hybrid genetic search algorithm.

\subsection{UAV for urban logistics}
\label{uav_intro}

The increasing adoption of UAVs is due to their autonomy and efficiency. They are also versatile and affordable in relatively large numbers, and these attributes are getting better with time. Recent advancements in drone technology exhibit improvements along multiple dimensions \citep{chen2021improved}. Upgrades in battery and propulsion technology allow enhanced manoeuvrability, longer flight times, greater altitude capabilities and higher payload capacity. For example, \cite{hecken2022structural} studies the construction of a heavy-lift UAV for low-altitude flight. The incorporation of machine learning and swarm technology permits greater degrees of autonomy and coordination in large groups. Progress in materials science and UAV design create scalability in drone manufacture, reduction in overall weight and better operational safety. \cite{kornatowski2020morphing} design a drone with collapsible rotor arms to improve safety while landing and aerodynamic efficiency during flight. Efforts are also being made to conceptualize novel UAV infrastructure, such as docking stations, precision landing sites, and innovative battery charging modes \citep{mohsan2023unmanned}. 
\\
\\
UAV deployment for traffic monitoring has witnessed growing attention from a research perspective. UAV patrols involve repeated deployments to specific parts of a city to monitor traffic flow based on volume \citep{chow2016dynamic}. Roadside surveillance by UAVs also extends to the conception of intelligent transport systems. Working as accident report agents, UAVs can alert nearby autonomous vehicles and assist rescue teams to minimize and mitigate the damage caused \citep{menouar2017uav}. Aerial data analysis can also reveal possible blockages at both traffic and railroad intersections \citep{congress2021identifying}. Visual and spatial information gathered by planned UAV flights can help make crucial decisions for disaster management. \cite{glock2020mission} explore route planning for UAVs to detect sampling locations and collect atmospheric samples after a fire or hazardous chemical accident. This data is then used to predict the atmospheric distribution of toxic substances in the proximity of the affected region.
\\
\\
There are many models for payload delivery using UAVs in the literature. UAVs can either be deployed independently, or in association with conventional vehicles such as trucks (with or without synchronization). This paper focuses on the former, i.e. the deployment of a homogeneous fleet of UAVs without ground vehicular support. However, for completeness, we delve into the relevant literature for both drone-only and truck and drone logistics, in that order. 

\subsubsection{Drone-only routing}

UAVs are cheaper to maintain, reduce labour and delay costs and have a lower carbon footprint compared to trucks. Amazon and Alibaba have been frontrunners in adopting UAVs for logistics, in moderate density urban regions and to deliver time-sensitive packages respectively. Many other companies across the world are following suit~\citep{han2022covid}. Additionally, a majority of customers in urban locations live close to a supermarket, and most delivery packages are relatively light-weight \citep{moshref2021applications}. This includes food, beverages, groceries and also medical and pharmaceutical supplies, all of which have a perishability constraint which makes UAVs an attractive choice.
\\
\\
\cite{levin2023branch} explore the drone delivery service planning problem and propose a novel branch-and-price algorithm for its solution. They also design a reservation heuristic to generate feasible solutions expediently. \cite{vlahovic2017implementing} perform a case study of a UAV delivery model for pharmaceutical supplies, exhibiting cost reductions and improvements in service responsiveness. Retail facility placements can be improved when network design is conducted under the availability of a UAV delivery service \citep{baloch2020strategic}. \cite{liu2019optimization} propose a mixed integer programming formulation for online meal delivery using UAVs. The objective function penalizes excess battery depletion, motivates fast delivery and minimizes unnecessary movement for re-routing. The model is executed in a progressive rolling fashion across discrete time-steps. \cite{chen2021improved} explore key strategic and tactical decisions for retailers implementing drone-based delivery operations, focusing on when to offer drone delivery and the appropriate pricing. Using a Markov decision process (MDP) framework, they introduce heuristic procedures to obtain near-optimal closed-form solutions efficiently. The aim is to assist online retailers in real-time decision-making regarding the feasibility and extent of offering drone-based delivery for specific product categories in various service zones. Additionally, they identify effective pricing strategies for drone-based deliveries. Payload delivery expands outside the realm of delivery to consumers. For example, UAVs can be used in precision agriculture, to simultaneously monitor crop sites and perform spraying tasks. This increases the efficiency of both pesticide and fertiliser use \citep{radoglou2020compilation}.
\\
\\
The study of drone-only logistics also attracts ML-enriched solution strategies. \cite{asadi2022markov} study drone delivery of medical supplies to hospitals, categorizing demand based on hospital distance from the drone hub. Each demand class includes a demand distribution and the required battery charge for successful delivery. The problem is formulated as a MDP where the system state is a vector containing the number of batteries at specific charge levels within the hub. The concise state expression allows reinforcement learning where a lookup table suffices for value function approximation. The action space involves charging batteries from lower to higher levels within decision epochs. This model optimizes battery usage, preserving charge for future deliveries and reducing replacement and charging costs. UAVs also have potential as a supplementary resource in two-echelon logistics networks. \cite{li2024parcel} study the last-mile delivery problem, where trucks transport goods from a primary depot to satellite depots, and UAVs complete the delivery to customers. Note that we categorize this study as drone-only logistics since trucks are treated as a stochastic externality. The stochastic nature of parcel arrival times at satellites is handled through Bayesian estimation, converting continuous arrival scenarios into discrete expected arrival intervals. A two-stage stochastic model is applied: parcel consolidation is first performed to release a subset for immediate delivery, and the remaining parcels are reserved for routing with future arrivals. Next, the route is treated as a particle in a particle swarm optimization problem. Reinforcement learning using search operators incrementally permutes the route in each iteration to discover better routes in the local neighbourhood. 

\subsubsection{Truck and drone logistics}

The integration of UAVs into the logistic service fleet (together with conventional vehicles) can result in considerable benefits \citep{rejeb2023drones}. Such a multi-modal delivery model decreases delivery times, overcomes traffic congestion and minimises human participation in the delivery process. \cite{dayarian2020same} investigate a same-day delivery model where a fleet of trucks is assisted by UAVs (through commodity resupply) at mutually determined meeting locations, allowing trucks to make longer trips before returning to the depot. \cite{gu2023dynamic} investigate a single truck-and-drone scenario for online demand, consisting of deterministic delivery requests along with online pickup requests. The drone revisits the truck for resupply when needed, such revisits dividing the route into segments. An insertion heuristic compares the costs of placing pickup customers at various segments of the initial route. A $15\%$ increase in profits is reported, with UAVs contributing a $50\%$ increase in the acceptance rate of online requests. \cite{yang2023planning} study a robust drone-truck delivery problem with customer time windows. They leverage the precision of drone routes to mitigate uncertainties in truck travel times caused by ground traffic congestion. Unlike \cite{gu2023dynamic}, here the truck remains at a customer location until the drone serves other customers and returns. Numerical experiments show that the drone-truck combination can safe-guard against the risk of failed deliveries due to unmet time windows. \cite{rave2023drone} broaden the truck-drone route planning framework by introducing deployable drone launching stations known as microdepots. In this setup, trucks leave a central distribution centre equipped with both onboard drones and microdepots. The authors address a truck-drone fleet planning problem for a single logistics service provider using adaptive large neighbourhood search. Experiments demonstrate that mixed fleets provide cost-savings over a truck-only fleet, and that the cost-optimal choice of fleet composition and delivery method involves considering all possible combinations.
\\
\\
Heterogeneous fleets logistics, such as truck and drone models, also welcome ML-supported solutions. \cite{chen2022deep} handle a mixed fleet of vehicles and UAVs to handle online requests. Network evolution is modelled as a MDP. Routing heuristics determine route feasibility by vehicle or drone, and deep reinforcement learning is then used to determine if a vehicle or a drone should be dispatched. \cite{ghiasvand2024data} tackle a two-echelon multi-trip truck-drone routing problem with uncertainty. Trucks make multiple trips from the depot to local distribution centres (LDCs), while UAVs make single trips from LDCs to customers. The research formulates a mixed-integer linear programming model to minimize total customer waiting times, managing uncertain customer demand through kernel-function based machine learning algorithms.
\\
\\
In an urban scenario, the operation environment of UAVs can mimic the road network. Buildings and other infrastructure can impede flight paths or require excessive energy to overcome aerially. But the airspace above roadways is generally free of obstructions. A possible drone flight network can be a map of the road traffic network itself, where links are above roads and nodes are above traffic intersections and several flight levels may be considered~\citep{levin2023branch}. We adopt this problem and context for the online DDSP.
    
\section{Problem description and formulation}
\label{prob_des}

We first introduce the deterministic demand formulation for the DDSP in Section~\ref{staticDDSP}. We next present the online demand model in Section~\ref{stoc_dem} and the MDP formulation in Section~\ref{MDP_form}.

\subsection{Deterministic DDSP}
\label{staticDDSP}

\begin{figure}[H]
    \centering
    \includegraphics[height=4cm]{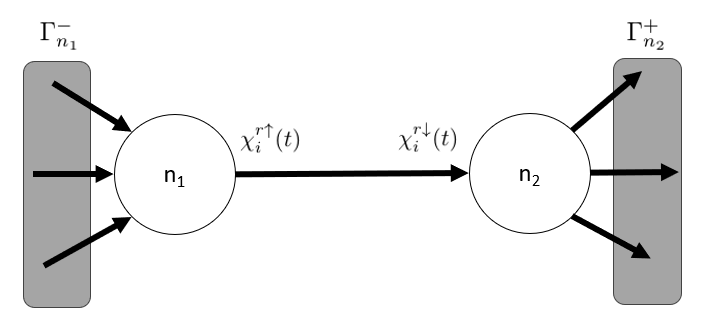}
    \caption{Link-Based Formulation}
    \label{fig:link_form}
\end{figure}

The deterministic (or static) version of the DDSP is examined by \cite{levin2023branch}  through an integer linear programming (ILP) formulation. This formulation is discussed below and serves as the foundation for the model and solution approach designed in this study for the online DDSP. 
\\
\\
Let $G = (N,A)$ be a network across which customers send transport requests. These requests, if accepted, are fulfilled by a fleet of drones. Every link $i \in A$ has a length $\rho_{i}$ and capacity $C_{i}$. All links are directed. The drones traverse the network through these links at a constant velocity $v$. Each node $n$ has a set of incoming and outgoing links $\Gamma^{-}_{n}$ and $\Gamma^{+}_{n}$ respectively. The time horizon is discretized into time steps $\{0,1,\ldots,t,\ldots,T\}$, where $T$ marks the end of the horizon. Customers submit requests to the network \textit{before} the time period of network operation $T$. Each customer request $r$  consists of (a) spatial parameters -- the request origin $o_{r}$ and destination $d_{r} $ such that $o_{r},d_{r} \in N$, (b) temporal parameters -- the earliest permissible time of departure $e_{r}$ and the arrival time window at the destination $[l_{r}, u_{r}] $ and (c) a profit $p_{r}$ earned upon fulfilling this request. The set of all requests that arrive over the entire time horizon is denoted by $\mathcal{R}$. In the static version, the network manager is aware of the complete request set $\mathcal{R}$ before the time horizon begins. 
\\
\\
\cite{levin2023branch} introduce a link-based network formulation. Every link is divided into an upstream and downstream segment (Figure \ref{fig:link_form}). The key decision variable $\chi^{r\uparrow}_{i}(t) \in \{0, 1\}$ indicates whether request $r$ occupies the upstream portion of link $i \in A$ at time $t$. The counterpart of $\chi^{r\uparrow}_{i}(t)$ for the downstream segment is $\chi^{r\downarrow}_{i}(t)$.

Two successive links connected by a node is called a $turn$. The decision variable $\gamma_{ij}^{r}(t) \in \{0,1\}$ indicates whether request $r$ occupies the turn from link $i$ to $j$ at time $t$. Note that the movement of the drone from the downstream segment of link $i$ to the upstream segment of link $j$ through the connecting node is assumed to be instantaneous. A single node can serve as the junction for many turns. For a given turn $(i,j)$, the set of all other turns sharing the same junction node is denoted by $\mathcal{C}_{ij}$. The decision variable $\phi_{ij}(t) \in \{0,1\}$ resolves turn conflicts at every node by indicating which turn is active at a given time $t$, consequently deactivating all other turns which share that node as a junction. 
\\
\\
\noindent Finally, the decision variable $z_{r} \in \{0,1\}$ indicates the acceptance or rejection of request $r$, and the objective is to maximize the total profit acquired. The ILP Formulation~\eqref{deterministicDDSP} is presented below:

\begin{subequations}
\begin{align}
\max \ Z = &\sum\limits_{r \in \mathcal{R}}z_{r} p_{r} && \\
\text{s.t.} &\sum\limits_{i \in \Gamma_{o_{r}}^{+}} \sum\limits_{t=e_{r}}^{T} \chi^{r\uparrow}_{i}(t) = z_{r} &&\forall r \in \mathcal{R} \label{eq:1b} \\
& \sum\limits_{i \in \Gamma_{d_{r}}^{-}} \sum\limits_{t=l_{r}}^{u_{r}} \chi^{r\downarrow}_{i}(t) = \sum\limits_{i \in \Gamma_{o_{r}}^{+}} \sum\limits_{t=e_{r}}^{T} \chi^{r\uparrow}_{i}(t)  &&\forall r \in \mathcal{R} \label{eq:1c} \\
&\chi^{r\downarrow}_{i}\left(t + \frac{\rho_i}{v}\right) = \chi^{r\uparrow}_{i}(t) &&\forall r \in \mathcal{R}, \forall i \in A, \forall t \in \{0,\ldots,T\} \label{eq:1d} \\
&\sum\limits_{r \in \mathcal{R}} \chi^{r\uparrow}_{i}(t) \leq C_{i} &&\forall i \in A, \forall t \in \{0,\ldots,T\} \label{eq:1e} \\
& \sum\limits_{r \in \mathcal{R}} \chi^{r\downarrow}_{i}(t) \leq C_{i} &&\forall i \in A, \forall t \in \{0,\ldots,T\} \label{eq:1f} \\
&\sum\limits_{i \in  \Gamma_{n}^{-}} \chi^{r\downarrow}_{i}(t) = \sum\limits_{i \in  \Gamma_{n}^{+}} \chi^{r\uparrow}_{i}(t) &&\forall n \in N, \forall r \in \mathcal{R}, \forall t \in \{0,\ldots,T\} \label{eq:1g} \\
&\gamma^{r}_{ij}(t) \leq \chi^{r\downarrow}_{i}(t) &&\forall r \in \mathcal{R},\forall (i,j) \in A^{2}, \forall t \in \{0,\ldots,T\} \label{eq:1h} \\
&\gamma^{r}_{ij}(t) \leq \chi^{r\uparrow}_{i}(t) &&\forall r \in \mathcal{R},\forall (i,j) \in A^{2}, \forall t \in \{0,\ldots,T\} \label{eq:1i} \\
& \gamma^{r}_{ij}(t) \geq \chi^{r\downarrow}_{i}(t) + \chi^{r\uparrow}_{i}(t) - 1 &&\forall r \in \mathcal{R},\forall (i,j) \in A^{2}, \forall t \in \{0,\ldots,T\} \label{eq:1j} \\
&\phi_{ij}(t) \geq \frac{1}{|\mathcal{R}|} \sum\limits_{r \in \mathcal{R}} \gamma^{r}_{ij}(t) && \forall (i,j) \in A^{2}, \forall t \in \{0,\ldots,T\} \label{eq:1k} \\
& \phi_{ij}(t) + \sum\limits_{(i',j') \in \mathcal{C}_{ij}}\phi_{i'j'}(t) \leq 1 &&\forall (i,j) \in A^{2}, \forall t \in \{0,\ldots,T\} \label{eq:1l} \\
&\gamma_{ij}(t) \in \{0,1\} &&\forall r \in \mathcal{R}, \forall (i,j) \in A^{2}, \forall t \in \{0,\ldots,T\} \label{eq:1m} \\
&\phi_{ij}^{r}(t) \in \{0,1\} && \forall (i,j) \in A^{2}, \forall t \in \{0,\ldots,T\} \label{eq:1n} \\
&\chi^{r\uparrow}_{i}(t), \chi^{r\uparrow}_{i}(t) \in \{0, 1\} &&\forall r \in \mathcal{R}, \forall i \in A, \forall t \in \{0,\ldots,T\} \label{eq:1o} \\
&z_{r} \in \{0, 1\} &&\forall r \in \mathcal{R} \label{eq:1p}
\end{align}
\label{deterministicDDSP}
\end{subequations}

Equations \eqref{eq:1b} and \eqref{eq:1c} ensure that only accepted requests are routed, and the route reaches the destination $d_{r}$ within the time window. Equation \eqref{eq:1d} monitors link travel time while equations \eqref{eq:1e}--\eqref{eq:1f} prevent capacity violations. Equation \eqref{eq:1g} conserves link flow, and equations \eqref{eq:1h}--\eqref{eq:1j} construct valid turns. Lastly, equations \eqref{eq:1k}--\eqref{eq:1l} resolve conflicts at turn junctions (i.e., nodes with non-empty $\Gamma_{n}^{+}$ and $\Gamma_{n}^{-}$) in the network.

\subsection{Online demand model}
\label{stoc_dem}

\begin{figure}[t]
    \centering
    \includegraphics[width=0.5\textwidth]{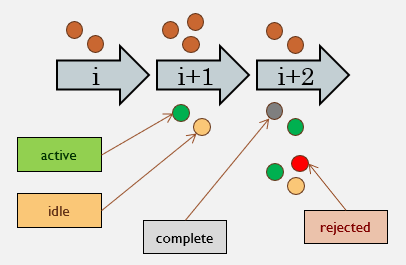}
    \caption{Online demand model: the circles represent requests while $i, i+1, i+2$ are time intervals}
    \label{fig:stochastic_demand}
\end{figure}

We build upon the deterministic DDSP to establish the online demand model. The time horizon is also divided into $I$ intervals, each of duration $D$, such that $DI = T$. The successive execution of all $I$ intervals concludes a single instance (or day) of network operation. The number of requests arriving at the network in each interval is governed by a Poisson process with rate $\lambda$. This implies that, on average, the expected total number of requests submit to the network during the entire time horizon is $\lambda I$. Within a specific interval, the time of submission of each request (by the customer) is assumed to be distributed uniformly.
\\
\\
Customers submit requests to the network during each interval $i \in I$. The current incoming requests in interval $i$ (denoted by $\mathcal{R}_{C}^{i} \subset \mathcal{R}$) can either be accepted or rejected. Those which are accepted are added to the set $\mathcal{R}_{\checkmark} \subset \mathcal{R}$ and must then be routed through the network. Once a request is accepted (assigned a drone and a planned route), it cannot later on be rejected or left unfulfilled. As soon as a drone initiates its route within the network, the corresponding request is termed $active$. Alternatively, if a request is accepted (and routed) but the assigned drone has not left the origin, the request is termed $idle$. Let $\mathcal{R}_{I} \subset \mathcal{R}_{\checkmark}$ and $\mathcal{R}_{A} \subset \mathcal{R}_{\checkmark}$ be the sets of idle and active requests respectively. All idle requests in the current interval are carried over to the next one, where their planned routes are subject to change. Once departed, no route changes are permitted and the drone is assumed to complete the delivery and drop off the network. Upon route completion the corresponding request is removed from  $\mathcal{R}_{A}$ and $\mathcal{R}_{\checkmark}$. The problem objective is to maximize the total profit earned by servicing requests in the given time horizon. This demand model is illustrated in Figure~\ref{fig:stochastic_demand}.

\subsection{Formulation as a Markov decision process}
\label{MDP_form}

The dynamic evolution of the network (arrival and acceptance of requests, followed by drone traffic routing and optimization) can be formulated as a Markov Decision Process (MDP). In an MDP, the problem is decomposed into an $environment$ and an $agent$. In the online DDSP, the network is the environment, and the network manager (who accepts and routes requests) is the agent. The environment produces $states$, based on which the agent takes $actions$, to receive $rewards$. The agent decides the best action at a given state. The environment subsequently controls the transition to the next state.

\subsubsection{State space}

The evolution of drone traffic in the network is represented as a sequential MDP moving through states. The state of the network captures all relevant information about its current structure. The state $s_{r}$ is constructed with the following components:
\begin{enumerate}
\item{$t_{r}$: The time step at which a new request arrives. This triggers the creation of a new state. Note that two successive requests may arrive several time steps apart i.e., $t_{r+1} - t_{r} \geq 1$.}
\item{$\vec{r}$: The request vector contains details of the request. These are: the request sequence number, origin node, destination node, earliest departure time, time window of arrival at the destination, and request profit. We write $\vec{r} = (r, o_{r}, d_{r}, e_{r}, [l_{r}, u_{r}], p_{r})$.}
\item{$\vec{\Theta}_{r}^{pre}$: The route vector $\vec{\theta}_{r}$ contains details of the drone route assigned to request $r$. The collection of all route vectors in state $s_{r}$ is $\vec{\Theta}_{r}^{pre}$. The composition of $\vec{\theta}_{r}$ reveals the request sequence number, the time of drone departure, and sequence of nodes in the route $(r, t_{o_{r}},  [o_{r},\ldots,d_{r}])$. For sake of brevity, the sequence of nodes in the route is also denoted by $\mu_{r} = [o_{r},\ldots,d_{r}]$. It can be known if the request $r$ is active or idle based on whether $t_{o_{r}}$ lies within the current interval or without. If a request becomes active, its route vector is immune to any further change. If a request is fulfilled by the time the new state is triggered, its corresponding $\vec{\theta}_{r}$ is removed from $\vec{\Theta}_{r}^{pre}$.}
\end{enumerate}

Hence the state $s_{r}$ is represented as the tuple $s_{r} = (t_{r}, \vec{r}, \vec{\Theta}_{r}^{pre})$. Every interval witnesses the creation of as many states as there are request arrivals in that interval. Between the acceptance (and subsequent routing) or rejection of a request and the arrival of the successive request, the only change in the network is the movement of active drones as per their routes. Since this is a deterministic change (the routes of active drones are known and fixed), a new state need not be declared until the arrival of the next request.

\subsubsection{Action space}

After the arrival of a new request $r$ triggers the creation of state $s_{r}$, the agent must take an action $a_{r}$ from the feasible action space $\mathcal{A}(s_{r})$. The action space allows the rejection of request $r$, but also contains all possible ways in which request $r$ can be routed through the network, provided four conditions are met:

\begin{enumerate}
\item{The drone must not exceed the capacity of a link pre-occupied by other active drones.}
\item{The drone must not, at any time-step, utilize a turn junction which is already utilized by an active drone at that time-step.}
\item{The drone must not depart the request origin before $e_{r}$.}
\item{The drone must reach the request destination within the time window $[l_{r}, u_{r}]$.}
\end{enumerate}

The action space $\mathcal{A}(s_{r})$ also permits re-routing idle requests from previous intervals in order to accommodate request $r$ or to achieve a better objective value. The action taken is then represented as the tuple $a_{r} = (z_{r}, \vec{\Theta}_{r}^{post})$ where $z_{r}$ is a binary decision regarding request $r$ and $\vec{\Theta}_{r}^{post}$ is the $updated$ collective route vector, which would include a route $\vec{\theta}_{r}$ for request $r$ if $z_{r} = 1$. At a given state, action selection is done by the agent's \textbf{policy}, which maps states to actions. Conventionally, a policy is represented as $\pi (a_{r}|s_{r})$, and governs the probability of taking action $a_{r}$ when in state $s_{r}$. 

\subsubsection{Reward and value functions}

Accepting a request increases the value of the objective function and should therefore indicate a reward. The reward received in state $s_{r} = (t_{r}, \vec{r}, \vec{\Theta}_{r}^{pre})$ with action $a_{r} = (z_{r}, \vec{\Theta}_{r}^{post})$ is simply:
\begin{equation}
R(s_{r}, a_{r}) = z_{r} p_{r}
\end{equation}

This gives a reward of $p_{r}$ for accepting a request and a reward of 0 for rejecting it. Consequently, the reward accumulated over the course of interval $i \in I$ is given by $\sum\limits_{r \in \mathcal{R}_{C}^{i}} R(s_{r}, a_{r}) = \sum\limits_{r \in \mathcal{R}_{C}^{i}} z_{r} p_{r}$. The overall cumulative reward is $\sum\limits_{t_{r} \in [0, T]} z_{r} p_{r}$, which is equivalent to $\sum\limits_{r \in \mathcal{R}} z_{r} p_{r}$. 
\\
\\
The state-action pair $(s_{r}, a_{r})$ captures the process of taking an action $a_{r}$ while in a state $s_{r}$. In the application of ML under the MDP framework, every state-action pair is assigned a relative value in comparison to all other pairs. This value indicates the expected final objective (or reward) that can be achieved given that the agent occupies state $s_{r}$ and takes action $a_{r}$. The learning model presented in this study seeks to estimate the value of taking a certain action while occupying a certain state. In other words, it seeks to discover and exploit the \textit{action value function $Q(s_{r}, a_{r})$}. Thus ``learning” is equivalent to improving the accuracy of the value estimations of $Q(s_{r}, a_{r})$ to a satisfactory degree.

\subsubsection{Online information and state transition}
 
The transition to the next state occurs under transition probabilities $P^{a_{r}}_{s_{r}\rightarrow s_{r}'}$. The transition probability determines which state the environment jumps to once the agent undertakes a specific action in the current state. There are many application scenarios where the transition probability matrix is either unknown to the agent, or the size of the state space prohibits its use. In such scenarios the agent can only control its policy, and has no prior knowledge of how the environment will react to its actions. The drone delivery network with online demand falls under this scenario. 
\\
\\
As mentioned, the density of request arrival across the time horizon is due to a Poisson process with rate $\lambda$. The origins and destinations of the requests are generated via spatial distributions, the time windows of arrival at the destinations are governed by a temporal distribution, and the request profits are sourced from a profit distribution. However, the agent or network manager has no prior information about any of these distributions or the Poisson process. This makes the problem online and the nature of request arrival becomes exogenous to the agent.  

\subsubsection{Objective function and optimal policy}

The objective is to maximize the cumulative profit earned from requests served over the entire horizon. This is identical to maximizing the cumulative reward acquired in the MDP formulation. In MP terminology, this becomes:

\begin{equation}
\text{max} \ Z = \sum\limits_{r \in \mathcal{R}}z_{r} p_{r}
\label{eq:obj}
\end{equation}

\noindent Given an (accurate) action value function $Q(s_{r}, a_{r})$, an optimal policy $\pi^{*}$ would aim to choose action $a_{r}$ in state $s_{r}$ such that the state-action pair $(s_{r}, a_{r})$ returns the maximum achievable $Q$-value in that state, i.e.:
\begin{equation}
\pi^{*} (a_{r}|s_{r}) = \begin{cases}
    1, & \text{if } a_{r} = \underset{a}{\arg\max} \ Q(s_{r},a) \\
    0, & \text{otherwise }
\end{cases}
\end{equation}

This is equivalent to stating that the policy, out of the set of all policies $\Pi$, that maximizes the objective function, i.e. the expected total reward gained over the course of the time horizon, is $\pi^{*}$:
\begin{equation}
\pi^{*} \in \underset{\pi \in \Pi}{\arg\max} \ \mathbf{E} \left[ \sum\limits_{r = 1}^{|\mathcal{R}|} R(s_{r}, A_{r}^{\pi}(s_{r})) \Big| s_{r-1}  \right] 
\label{eq:pi_star}    
\end{equation}

Here $A_{r}^{\pi}(s_{r})$ is the action taken by policy $\pi$ in state $s_{r}$. The state $s_{0}$ depicts the initial empty network where no request has yet arrived. Note that the limits of summation in \eqref{eq:obj} and \eqref{eq:pi_star} are slightly different because, unlike in \eqref{eq:obj}, in \eqref{eq:pi_star} we emphasize operating upon requests sequentially to demonstrate the evolution of the MDP.

\section{Solution methodology}
\label{sol_met}

A rolling implementation of Formulation~\eqref{deterministicDDSP} in consecutive intervals $i \in I$ is developed to solve the online DDSP. Section \ref{rolling_milp} first outlines the rolling ILP implementation as a solution framework. This is followed by Section \ref{sof} which illustrates the motivation to construct a surrogate ILP objective function. The determination of predictive prescriptions via a training scheme is provided in Section~\ref{training}. The resulting data-driven optimization approach and the proposed policy for solving the online DDSP are summarized in Section~\ref{policy}.

\subsection{Myopic ILP}
\label{rolling_milp}

Five modifications need to be made to Formulation~\eqref{deterministicDDSP} for the dynamic online DDSP. First, since Formulation~\eqref{deterministicDDSP} was developed for the static DDSP, all time-indexed variables and constraints extend over the entire time horizon $T$. However, in the online demand context, the ILP will be solved repeatedly for each interval $i \in I$ in a rolling fashion. This requires that, in a given interval, the time-indexed variables and constraints extend from the beginning to the end of that interval. Hence the limits of $t$ change from $\{0,\ldots,T\}$ to $\{(i-1)D,\ldots,iD\}$ for interval $i$.
\\
\\
Second, since idle requests are carried over from previous intervals, any change in their $z_{r}$ values is not permissible. Requests once accepted must be fulfilled. The same is true for active requests from preceding intervals as well. This implies that in each interval,
\begin{subequations}
\begin{equation}
z_{r} = 1 \qquad \forall r \in \mathcal{R}_{I} \cup \mathcal{R}_{A}
\label{eq:6a}
\end{equation}
\\
Third, the objective function needs to be implemented for all requests which have arrived in the current interval $i$, and for idle requests from previous intervals. The latter is because idle requests are still a part of the solution process, since their routes are subject to re-optimization. This is captured as:
\begin{equation}
\text{max} \ Z = \sum\limits_{r \in \mathcal{R}_{O}^{i}} z_{r} p_{r} \quad \text{where} \ \mathcal{R}_{O}^{i} = \mathcal{R}_{C}^{i} \cup \mathcal{R}_{I} \ \forall i \in I
\label{eq:6b}
\end{equation}

Fourth, since idle requests may still be re-routed, their earliest time of departure $e_{r}$ cannot precede the start of the current interval $i$. An idle request cannot be re-routed such that its new route begins at a time-step that has already passed.
\begin{equation}
e_{r} =
\begin{cases}
  (i-1)D, \quad \text{if} \ e_{r} \leq (i-1)D \\
  e_{r}, \quad \text{otherwise} \\  
\end{cases} \forall r \in \mathcal{R}_{I}
\label{eq:6c}
\end{equation}

Finally, the unfinished portion of active request routes which will be undertaken in the current interval must be preserved. Suppose we are in interval $i$ and the last updated route vector for (currently) active request $r$ is $\vec{\theta}_{r} = (r, t_{o_{r}},  [o_{r},\ldots,d_{r}])$. Let the number of nodes in its path $\mu_{r}$ be $P$ and let $\omega(m,n)$ denote the directed link from node $m$ to node $n$. Then the definitive path for request $r$ is $\mu_{r} = [o_{r}, n_{2},\ldots,n_{k},\ldots,n_{P-1},d_{r}]$ where $n_{1} = o_{r}$ and $n_{P} = d_{r}$. The following impositions must be made for all such active requests $r \in \mathcal{R}_{A}$:
\begin{equation}
\begin{rcases}
  \chi^{r\downarrow}_{\omega(n_{k-1},n_{k})}(t_{n_{k}}^{r}) = 1, \quad \text{if} \ n_{k} \neq o_{r} \\
  \chi^{r\uparrow}_{\omega(n_{k},n_{k+1})}(t_{n_{k}}^{r}) = 1, \quad \text{if} \ n_{k} \neq d_{r}\\
  \gamma_{\omega(n_{k-1},n_{k})\omega(n_{k},n_{k+1})}(t_{n_{k}}^{r}) = 1, \quad \text{if} \ n_{k} \notin \{o_{r}, d_{r}\}\\
\end{rcases} \forall \ n_{k} \in \mu_{r} \mid (i-1)D \leq t_{n_{k}}^{r} \leq iD
\label{eq:6d}
\end{equation}
\end{subequations}

$t_{n_{k}}^{r}$ denotes the time at which the drone responsible for request $r$ reaches node $n_{k}$. $t_{n_{k}}^{r}$ can be determined since $t_{o_{r}}$, link lengths and drone velocity are known. Solving Formulation~\eqref{deterministicDDSP} with these five modifications at every interval completes the rolling ILP implementation. The pseudo-code is presented in Algorithm~\ref{alg:roll-milp}. Since this model does not consider the impact of current decisions on the fate of requests submitted in future intervals, and only focuses on the current interval, we hereby refer to it as the Myopic ILP.

\begin{algorithm}[]
\caption{Pseudo-code for Myopic ILP}
\label{alg:roll-milp}

\For{$d \leftarrow 1$ \KwTo $N_{days}$}{
  initialize request set $\mathcal{R}_{\checkmark}$ and collective route vector $\vec{\Theta}$ \\
  \For{$i \leftarrow 1$ \KwTo $I$}{      
     initialize request sets $\mathcal{R}_{A}, \mathcal{R}_{I}, \mathcal{R}_{C}^{i}$ \\  
     \textcolor{magenta}{\tcp{request segregation}}
    \For{$r \in \mathcal{R}_{\checkmark}$}{
      \eIf{$t_{o_{r}} \leq (i-1)D$}
            {$\mathcal{R}_{A} \leftarrow \mathcal{R}_{A} \cup \{r\}$}
            {\If {$e_{r} \leq (i-1)D$}{
     $e_{r} \leftarrow (i-1)D$ 
     }
      $\mathcal{R}_{I} \leftarrow \mathcal{R}_{I} \cup \{r\}$}      
    }
    populate $\mathcal{R}_{C}^{i}$ \\
    \textcolor{magenta}{\tcp{model formulation}}
    add constraints \eqref{eq:1b} -- \eqref{eq:1p} with $t \in \{(i-1)D,\ldots,iD\}$ wherever applicable \\
    \For{$r \in \mathcal{R}_{I} \cup \mathcal{R}_{A}$}{
    add constraint $z_{r} = 1$
    }
    \For{$r \in \mathcal{R}_{A}$}{
    add constraints \\
    $\begin{rcases}
  \chi^{r\downarrow}_{\omega(n_{k-1},n_{k})}(t_{n_{k}}^{r}) = 1, \ \text{if} \ n_{k} \neq o_{r} \\
  \chi^{r\uparrow}_{\omega(n_{k},n_{k+1})}(t_{n_{k}}^{r}) = 1, \ \text{if} \ n_{k} \neq d_{r}\\
  \gamma_{\omega(n_{k-1},n_{k})\omega(n_{k},n_{k+1})}(t_{n_{k}}^{r}) = 1, \ \text{if} \ n_{k} \notin \{o_{r}, d_{r}\}\\
\end{rcases} \forall \ n_{k} \in \mu_{r} \mid (i-1)D \leq t_{n_{k}}^{r} \leq iD$
    }   
\vspace{5pt}
   add objective $\max \ Z = \sum\limits_{r \in \mathcal{R}_{O}^{i}} z_{r} p_{r}$ \\
   solve modified Formulation~\eqref{deterministicDDSP}\\
   \textcolor{magenta}{\tcp{update/add drone routes}}
  \For{$r \in \mathcal{R}_{O}^{i}$}{
  \If {$r \in \mathcal{R}_{I}$}{
     update $\vec{\theta_{r}} \leftarrow \ (r, t_{o_{r}}, \mu_{r})$ 
     }
  \If {$r \in \mathcal{R}_{C}^{i}$}{
  \If {$z_{r}^{sol} = 1$}{
    $\mathcal{R}_{\checkmark} \leftarrow \mathcal{R}_{\checkmark} \cup \{r\}$ \\
    $\vec{\theta_{r}} \leftarrow \ (r, t_{o_{r}}, \mu_{r})$ \\
    $\vec{\Theta} \leftarrow \vec{\Theta} \cup \{\theta_{r}\}$ 
  }
  }
  }
  \textcolor{magenta}{\tcp{remove completed requests}}
  \For{$r \in \mathcal{R}_{\checkmark}$}{
  \If {$t_{d_{r}} \leq iD$}{
  $\mathcal{R}_{\checkmark} \leftarrow \mathcal{R}_{\checkmark} \backslash \{r\}$ \\
  $\vec{\Theta} \leftarrow \vec{\Theta} \backslash \{\theta_{r}\}$
  } 
  }
}
}
\end{algorithm}

\subsection{Surrogate objective function}
\label{sof}

Maximizing profit acquisition at each time interval ignores future demand. Taking future demand into account requires making decisions that are sub-optimal in the current time interval but increase the cumulative reward by contributing to the objective of a later interval. We tolerate a sub-optimal reward per interval $\sum\limits_{r \in \mathcal{R}_{C}^{i}}R(s_{r}, a_{r})$ to achieve a greater reward-to-go $\sum\limits_{r \in \mathcal{R}}R(s_{r}, a_{r})$. This sub-optimality comes from allocating some resources to the future, thus preventing their use in the present. In other words, a sacrifice is made in the present to reap greater rewards in the future. 
\\
\\
We define available link capacity per unit time as the network resource. This is captured as $C_{l} - \sum\limits_{r \in \mathcal{R}_{O}^{i}} \chi^{r\uparrow}_{l}(t)$, where the capacity of link occupied by idle and current requests is deducted from its total capacity. The strategy then translates into (intentionally) creating reserve capacity in the network at the current interval, which prevents that capacity from being utilized at the current interval. This may lead to rejection of requests or longer routes in the current interval, but is expected to yield greater profit accumulation in future intervals. We motivate the choice of using link capacity as a resource as follows: first, it is a quantitative definition, and link capacity at a given time can be easily measured during the execution of the optimization program; second, the spatial distribution of origins and destinations of customer requests, combined with conflict resolution of drone routes, leads to a heterogeneous resource space. In other words, the available capacities on different links have different value to the network. This creates a rich decision set, providing a sizeable margin with which to outperform myopic policies. 
\\
\\
The priority of available capacity on link $l \in A$ at time interval $i$ is denoted $\beta_{i,l}$. This indicates the significance of the link in the network. If a link has higher priority, then any reserve capacity on it is relatively more important to network operation and route creation, compared to other links. The formulation of a data-driven prescription of $\beta_{i,l} $ is described in Section~\ref{ddpp}. Reserve capacity is generated on network links in proportion to their priority, since slack created on higher priority links will have greater impact on maximizing profit in future intervals. 
\\
\\
Indefinite and uncontrolled slack creation will adversely affect network performance and reduce cumulative profit. The reserve capacity being created in the present should at some point be capitalized upon. Hence the balance between profit maximization and ($\beta_{i,l}$ proportionate) slack creation must be controlled. We construct a new objective function to incorporate the sacrificial strategy. Since this objective function is intended to replace the myopic objective function of the Myopic ILP, it is called the surrogate objective function and comprises two components:
\begin{equation}
\max \ Z = \sum\limits_{r \in \mathcal{R}_{O}^{i}} z_{r} p_{r} \ + \ \alpha_{i} \sum\limits_{l \in A} \beta_{i, l} \sum\limits_{t = (i-1)D}^{iD} \left[C_{l} - \sum\limits_{r \in \mathcal{R}_{O}^{i}} \chi^{r\uparrow}_{l}(t)\right]
\label{eq:7}  
\end{equation} 

The DDSP formulation with the surrogate objective function is named the Surrogate ILP. The first component is profit maximization, and the second is slack creation based on link priority. For every link $l \in A$, we first compute the slack capacity on $l$ during the current time interval, and then multiply it with the relative priority of that link, $\beta_{i,l}$. We also define a parameter $\alpha_{i}$ which influences the balance between the two components. The value of $\alpha_{i}$ chosen for this interval determines to what extent reserve capacity will be created on network links (relative to profit maximization). Note that independent of the value of $\alpha_{i}$, the slack generated on a link will be proportional to the relative priority of that link. Unless, of course, if $\alpha_{i} = 0$, in which we case the surrogate objective is reduced to the myopic objective of pure profit maximization and $\beta_{i,l}$ has no role left to play. 
\\
\\
In the computational experiments, both components of the surrogate objective function are standardized, so that $\alpha_{i}$ can effectively control the balance. The first component -- which is essentially the myopic objective function 
-- is standardized by division with the expected total profit over the entire time horizon $p_{avg}\mathbf{E}[|\mathcal{R}|]$. This component is hereafter referred to as the Profit Component. The second component is standardized by division with the total available network capacity in that interval $(\sum\limits_{l \in A} C_{l})D$. This component is hereafter referred to as the Slack Component.  

\subsection{Training procedure of relative link priority ($\beta_{i,l}$)}
\label{training}

We next explain the training procedure used to determine predictive prescriptions of the relative link priority ($\beta_{i,l}$) which is used in the Surrogate ILP. The training procedure is articulated across two main steps: i) the synthesis of training data via lookahead simulation (Section~\ref{ddpp}) and ii) the training of a supervised ML model for predicting relative link priority (Section~\ref{train_scheme}). To scale up the training process, we use heuristic algorithm to solve the Myopic ILPs that arise therein (Section~\ref{reserve}).

\subsubsection{Training data synthesis via lookahead simulation}
\label{ddpp}

We use lookahead simulations to synthesize training data that will be used to train a predictive model which itself will be queried during the execution of the  Surrogate ILP policy.
\\
\\
The priority of a link can be correlated to the frequency of link traversal by drones in a suitable span of time. Given a duration of network operation (e.g. $\delta$ intervals), crucial high priority links will experience large throughput and are likely to create bottlenecks which impede request acceptance. The link priority $\beta_{i,l}$ is therefore formally defined as the number of times link $l$ is traversed by incoming requests between the beginning of interval $i$ and the end of interval $i+\delta$.
\begin{equation}
\beta_{i,l} = \sum\limits_{q = i}^{i+\delta} \ \sum\limits_{r \in \mathcal{R}_{C}^{q}} \ \sum\limits_{t = (q-1)D}^{qD} \chi^{r \uparrow, sol}_{l}(t)
\label{eq:8}
\end{equation}

Here $\chi^{r \uparrow, sol}_{l}(t)$ represents the routing solution value of the upstream-link variable $\chi^{r \uparrow}_{l}(t)$ for link $l$ and time $t$ at the end of each interval. We also define a vector $\bm{\vec{B}}$ to store link priority values for all intervals. 

\begin{equation}
\bm{\vec{B}} = \begin{bmatrix}
    \bm{\vec{\beta}_{i}}
\end{bmatrix}^{i \in \{1, \ldots, I\}} = \begin{bmatrix}
    \beta_{i,l}
\end{bmatrix}_{l \in A}^{i \in \{1, \ldots, I\}}
\label{eq:9}
\end{equation}

Link priority values as defined above are a function of drone traffic evolution in the network, which in turn depends on customer demand. From a snapshot of the network state at a certain time-step, the short-term future of the network may be estimated, provided historical patterns of customer demand are known. The introduction of the ML framework serves this purpose. We develop a supervised ML model which takes as input a snapshot of the network (number of drones occupying each link) and returns as output an estimate of link priorities for the near future. This input snapshot of link occupation is captured for the beginning of every time interval and is denoted by:
\begin{equation}
s_{i, l} = \sum\limits_{r \in \mathcal{R}_{A}} \ \sum\limits_{t = (i-1)D - \Delta}^{(i-1)D + \Delta} \ \chi^{r \uparrow, sol}_{l}(t) \ \mathds{1}[t \leq (i-1)D < t + \Delta]  
\label{eq:10}
\end{equation}

where $\Delta = \frac{\rho_{l}}{v}$ is the time taken by a drone to traverse link $l$. Eq.~\eqref{eq:10} increments the value of $s_{i,l}$ for every active drone occupying link $l$ at the beginning of interval $i$ (i.e., $t = (i-1)D$). We define a collection vector to store link occupancy values for all intervals.
\begin{equation}
\bm{\vec{S}} = \begin{bmatrix}
    \bm{\vec{s}_{i}}
\end{bmatrix}^{i \in \{1, \ldots, I\}} = \begin{bmatrix}
    s_{i,l}
\end{bmatrix}_{l \in A}^{i \in \{1, \ldots, I\}} 
\label{eq:11}
\end{equation}

$\bm{\vec{S}}$ serves as the input feature vector dataset and $\bm{\vec{B}}$ as the target output vector dataset while training the ML model. Once trained, it is deployed in the Surrogate ILP to estimate the values of $\beta_{i,l}$ which are subsequently used in the surrogate objective function~\eqref{eq:7}. Note that in the training phase, all $\beta_{i,l}$ -- and hence $\bm{\vec{B}}$ -- are calculated via Eq.~\eqref{eq:8} to produce the output targets of the training dataset. In the testing and execution of the Surrogate ILP, however, the values of $\beta_{i,l}$ are predicted by the ML model. 
\\
\\
In the training scheme, the number of training intervals is designed to be much larger than that in a test instance ($I_{training} >> I$). This is equivalent to processing several test instances successively without emptying the network between two consecutive test instances. At the beginning of each training interval $i$, the link occupancy on each link $l$ by currently active requests is calculated and stored in $s_{i,l}$.  For this we define the function \textsc{LinkActivity} which computes $s_{i,l}$ as per Eq.~\eqref{eq:10}. This function takes two inputs -- the route $\mu_{r}$ of an active request and the link occupancy vector $\vec{\bm{s_{i}}}$ for the current interval. It iterates through the route to determine which link the UAV occupied at the beginning of the interval ($t = (i-1)D$) and increments the corresponding $s_{i,l}$ in $\vec{\bm{s_{i}}}$. Next, a lookahead is performed by virtually simulating the future of the network for the next $I_{virtual}$ number of intervals. For these intervals, virtual requests $\mathcal{R}_{V}$ are created, imitating historical data. These virtual requests are then routed through the network, obeying link capacity and conflict resolution constraints in the presence of active requests. The routes of successfully serviced virtual requests then contribute to the values of $\beta_{i,l}$ for interval $i$. For a link $l$, $\beta_{i,l}$ is equal to the number of times link $l$ was traversed by virtual requests in the lookahead simulation. This captures the perceived priority of link $l$ in the foreseeable near future, from the perspective of the network at interval $i$, in accordance with Eq.~\eqref{eq:8}, with $\delta = I_{virtual}$. We define the function \textsc{BetaUpdate}, which takes two inputs: the route $\mu_{r}$ of a virtual request and the link priority vector $\vec{\bm{\beta_{i}}}$ for the current interval, and updates $\beta_{i,l}$ for all the links that this request traverses according to $\mu_{r}$ . 
\\
\\
Once the virtual lookahead is complete, all virtual request paths are erased from the network. Thereafter, the process resembles the Myopic ILP. Idle requests from previous intervals and incoming requests of current interval $i$ are routed. Successfully accepted requests are stored in $\mathcal{R}_{\checkmark}$. Successful paths $\mu_{r}$ for a given request $r$ are also stored in a route vector $\theta_{r} = (r, t_{o_{r}}, \mu_{r})$ and added to the collective route vector $\vec{\Theta}$. Then the next interval $i+1$ is initiated.

\subsubsection{$k$-Nearest Neighbors multi-output regressor}
\label{train_scheme}

We use $k$-Nearest Neighbors ($k$NN) regression to train a predictive model for estimating relative link priority based on network snapshots. The $k$NN regressor is a supervised ML algorithm used for prediction tasks. It predicts the output target value of a new data point by considering the average or weighted average of its $k$ nearest neighbors in the input feature space. The $k$NN algorithm utilizes a defined distance metric, such as the Euclidean distance, to measure the similarity between data points in the input feature space. By identifying the $k$ training feature vectors closest to a new data point, the algorithm establishes the nearest neighbors of the given data point. The prediction of the target value of the new data point is made by averaging or weighting the target values associated with its $k$ nearest neighbors, leveraging the assumption that closer points in the feature space exhibit similar target values. The $k$NN multi-output regressor extends the concept of the $k$NN regressor to handle multiple target values for each data point. Notably, the $k$NN regressor is non-parametric and makes minimal assumptions about the underlying data distribution, making it more versatile in handling diverse datasets. It can also capture non-linear relationships between input features and output targets. We use the training data synthesized using the lookahead simulation procedure described in Section~\ref{ddpp} to train the $k$NN regressor. 

\subsubsection{Reservation heuristic}
\label{reserve}

\begin{figure}[t]
    \centering
    \includegraphics[height=7cm]{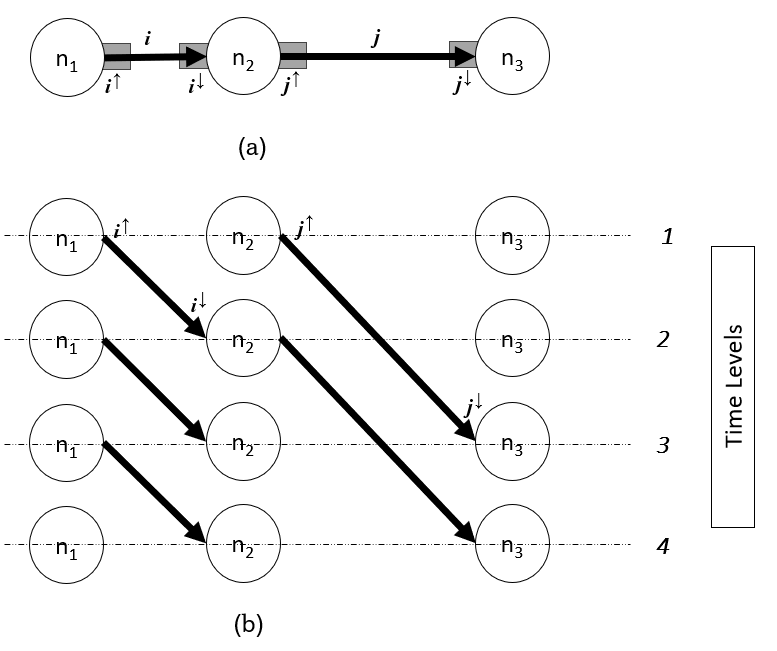}
    \caption{(a) Original Network (b) Time-expanded Network}
    \label{fig:time_expanded_network}
\end{figure}

In practice, the training phase has an associated time budget, implying that link priority values must be estimated expediently. Solving the Myopic ILP for $I_{training}$ intervals demands excessive time and memory. \cite{levin2023branch} propose a reservation heuristic which routes incoming customer requests quickly. This heuristic is used within the training phase. We illustrate the reservation heuristic and formalize the training scheme in Algorithm~\ref{alg:beta-train}.
\\
\\
The original spatial network $G = (N,A)$ can be converted into a spatio-temporal network $G^{T}$ as demonstrated in Figure~\ref{fig:time_expanded_network}. $G^{T}$ is the time-expanded dual graph of $G$. To construct $G^{T}$, the nodes and links are extended across time through discrete time levels and connections are made on the basis of link travel times. Suppose a link in $G$ connects node $n_{1}$ to node $n_{2}$ and can be traversed in $\Delta t$ time steps. Then it connects two nodes $\Delta t$ time levels apart in $G^{T}$. The graph $G^{T}$ can be viewed as a network of node-time tuples ($n,t$), where $n$ denotes the node and $t$ is the time level on which it is located. The reservation heuristic works by routing customer requests on $G^{T}$ using a shortest-path finding algorithm. The route is then ``reserved'' for that request by securing link capacity in space-time, making that capacity unavailable for other requests. Since turn conflict constraints must also be respected in $G^{T}$, the heuristic removes connectivity between node-time tuples whenever constraints \eqref{eq:1e}--\eqref{eq:1f} and \eqref{eq:1h}--\eqref{eq:1l} are violated. Algorithm~\ref{res-h} defines the function \textsc{NewReservation} which discovers a path $\mu_{r}^{T}$ for request $r$, makes the relevant reservations in $G^{T}$ and deduces the corresponding path $\mu_{r}$ in $G$. Two additional functions are also required in the training phase. A request may already have a path, but reservations for this path may either need to be made or erased. We use functions \textsc{Reserve}($\mu_{r}^{T}, G^{T}$) and \textsc{EmptyReservation}($\mu_{r}^{T}, G^{T}$) for these purposes respectively. They follow a similar procedure to \textsc{NewReservation} and are hence not illustrated.

\begin{algorithm}[H]
\KwIn{$r, G_{T}$}    
\KwOut{$\mu_{r}$}
    \SetKwFunction{FMain}{\textsc{NewReservation}}
    \SetKwProg{Fn}{Function}{:}{}
    \Fn{\FMain{$ r, G_{T} $}}{    
        find space-time path $\mu_{r}^{T}$ in $G_{T}$ satisfying \\
        \textbf{1}: previous reservations for \eqref{eq:1e}--\eqref{eq:1f} and \eqref{eq:1h}--\eqref{eq:1l} \\
        \textbf{2}: spatial and temporal constraints of request $r$ \\
        \eIf{$\mu_{r}^{T}$ \textsc{found}}{
            make reservations (break connectivity) in $G_{T}$ as per \eqref{eq:1e}--\eqref{eq:1f} and \eqref{eq:1h}--\eqref{eq:1l} for $\mu_{r}^{T}$ \\
            create $ \mu_{r} $ in $G$ from $\mu_{r}^{T}$ \\
            \textbf{return} $\mu_{r}$
        }{
            \textbf{return} \textsc{false}
        }    
    }
\textbf{End Function}
\caption{Pseudo-code for reservation heuristic}
\label{res-h}
\end{algorithm}

\begin{algorithm}[H]
\caption{Training data synthesis for predicting relative link priority}
\label{alg:beta-train}
initialize request sets $\mathcal{R}_{\checkmark}$ and collective route vector $\vec{\Theta}$ \\
set $s_{i,l} \leftarrow 0$ and $\beta_{i,l} \leftarrow 0 \ \forall \ l \in A, \forall \ i \in \{1,\ldots,I_{training} \}$ \\
create $G^{T}$ from $G$ \\
  \For{$i \leftarrow 1$ \KwTo $I_{training}$}{ 
  initialize request sets $\mathcal{R}_{A}, \mathcal{R}_{I}, \mathcal{R}_{C}^{i}$ \\ 
    \textcolor{magenta}{\tcp{request segregation as in Algorithm~\ref{alg:roll-milp} (lines 6 to 12)}} 
    \textcolor{magenta}{\tcp{reserve active requests and record link traversals}} 
    \For{$r \in \mathcal{R}_{A}$}{
     \textsc{Reserve}($\mu_{r}^{T}, G_{T}$) \\
     \textsc{LinkActivity}($\mu_{r}, \vec{\bm{s_{i}}}$) \\     
    }
    \textcolor{magenta}{\tcp{virtual requests}}   
    \For{$i_{v} \leftarrow 1$ \KwTo $I_{virtual}$}{    
    create virtual requests for interval $i_{v}$ and store in $\mathcal{R}_{V}$\\    
    }    
    \For{$r \in \mathcal{R}_{V}$}{     
    \If {\textnormal{\textsc{NewReservation(}}$r, G_{T} $\textnormal{\textsc{)}}}{    
     \textsc{BetaUpdate}($\mu_{r}, \vec{\bm{\beta_{i}}}$) \\ 
     }
    }
    \For{$r \in \mathcal{R}_{V}$}{
     \textsc{EmptyReservation}($\mu_{r}^{T}, G_{T} $)\\
     }
    \textcolor{magenta}{\tcp{solve real requests of current interval and update/add drone routes}}
    populate $\mathcal{R}_{C}^{i}$ \\
    \For{$r \in \mathcal{R}_{O}^{i}$}{
    $\mu_{r} \leftarrow$ \textsc{NewReservation(}$r, G_{T} $\textsc{)} \\
    \If {$r \in \mathcal{R}_{I}$}{
    update $\theta_{r} \leftarrow \ (r, t_{o_{r}}, \mu_{r})$
    }
    \If {$r \in \mathcal{R}_{C}^{i}$}{
    \If {$ \exists \mu_{r}$}{
    $\mathcal{R}_{\checkmark} \leftarrow \mathcal{R}_{\checkmark} \cup \{r\}$ \\
     $\theta_{r} \leftarrow \ (r, t_{o_{r}}, \mu_{r})$ \\
    $\vec{\Theta} \leftarrow \vec{\Theta} \cup \{\theta_{r}\}$\\
     }
    }  
    }  
     \textcolor{magenta}{\tcp{empty all reservations to prepare for next interval}}
    \For{$r \in \mathcal{R}_{O}^{i} \cup \mathcal{R}_{A}$}{
     \textsc{EmptyReservation}($\mu_{r}^{T}, G_{T}$)\\
     }
  \textcolor{magenta}{\tcp{remove completed requests as in Algorithm~\ref{alg:roll-milp} (lines 33 to 36)}}
}
\textbf{return} $\bm{\vec{B}}$ and $\bm{\vec{S}}$
\end{algorithm}

\begin{algorithm}[H]
\caption{Surrogate ILP policy}
\label{alg:surr-milp}
\For{$d \leftarrow 1$ \KwTo $N_{days}$}{
  initialize request set $\mathcal{R}_{\checkmark}$ and collective route vector $\vec{\Theta}$ \\
  initialize $\alpha$-profile $\begin{bmatrix}
    \alpha_{i}
\end{bmatrix}^{i \in \{1,\ldots,I\}}$\\
  set $s_{i,l} \leftarrow 0$ and $\beta_{i,l} \leftarrow 0 \ \forall \ l \in A, \forall \ i \in \{1,\ldots,I\}$ \\
  \For{$i \leftarrow 1$ \KwTo $I$}{      
     initialize request sets $\mathcal{R}_{A}, \mathcal{R}_{I}, \mathcal{R}_{C}^{i}$ \\  
     \textcolor{magenta}{\tcp{request segregation as in Algorithm~\ref{alg:roll-milp} (lines 6 to 12)}}
    \For{$r \in \mathcal{R}_{A}$}{     
     \textsc{LinkActivity}($\mu_{r}, \vec{\bm{s_{i}}}$) \\     
    }        
     \textcolor{magenta}{\tcp{query ML component for link priorities \& standardize}}
    $\bm{\vec{\beta_{i}}} \leftarrow$ \textsc{$k$NN}($\bm{\vec{s_{i}}}$) \\
    $\beta^{\max} \leftarrow$ $\max(\bm{\vec{\beta_{i}}})$ \\
    $\bm{\vec{\beta_{i}}} \leftarrow  \begin{bmatrix}
    \frac{\beta_{i, l}}{\beta^{\max}}
\end{bmatrix}_{l \in A}$\\
    populate $\mathcal{R}_{C}^{i}$ \\
    \textcolor{magenta}{\tcp{model formulation as in Algorithm~\ref{alg:roll-milp} (lines 15 to 20)}}
\vspace{-5pt}
   add surrogate objective $\max \ Z = \sum\limits_{r \in \mathcal{R}_{I}}z_{r} p_{r} \ + \alpha_{i} \sum\limits_{l \in A} \beta_{i, l} \sum\limits_{t = (i-1)D}^{iD} \left[C_{l} - \sum\limits_{r \in \mathcal{R}_{I}} \chi^{r\uparrow}_{l}(t)\right]$ \\
\vspace{5pt}
   solve modified Formulation~\eqref{deterministicDDSP} with surrogate objective\\
   \textcolor{magenta}{\tcp{update/add drone routes as in Algorithm~\ref{alg:roll-milp} (lines 24 to 31)}}
  \textcolor{magenta}{\tcp{remove completed requests as in Algorithm~\ref{alg:roll-milp} (lines 33 to 36)}}
}
}
\end{algorithm}

\subsection{Surrogate ILP policy} 
\label{policy}

The Surrogate ILP policy is executed after completing the training procedure described in Section~\ref{training}. The Surrogate ILP policy iteratively solve the Surrogate ILP where the surrogate objective function is parameterized by the predictive prescriptions for relative link priority ($\beta_{i,l}$) and the balance parameter ($\alpha_{i}$). Recall that the ratio of profit maximization to link priority-proportional slack creation is controlled via the balance parameter $\alpha_{i}$. The magnitude of this balance parameter may be varied from one time interval to the next. Since $\alpha_{i}$ is coupled with the slack component, a greater magnitude favours slack creation and a smaller one favours profit maximization. Too large a value of $\alpha_{i}$, sustained over many intervals, can impose a vacuum in the network which degrades cumulative profit acquired by injecting excessive reserve capacities. Too small a value can make its contribution insignificant, and a negative value is likely to suffocate the network by overloading link capacities. It follows that there exists a range of desirable values for $\alpha_{i}$, which we explore in our numerical experiments. The variation (or lack thereof) in the values of $\alpha_{i}$, denoted by the sequence $\begin{bmatrix}
    \alpha_{i}
\end{bmatrix}^{i \in \{1, \ldots, I\}}$, is called the \textit{$\alpha$-profile}. In our numerical experiments, three types of $\alpha$-profiles are investigated -- constant, step-wise and continuous. 
\\
\\
The pseudo-code the Surrogate ILP policy is presented in Algorithm~\ref{alg:surr-milp}. At the start of every interval, the current active requests (which were accepted in previous intervals) are segregated. The function \textsc{LinkActivity} populates the link occupancy vector $\vec{\bm{s_{i}}}$ in accordance with the positions of active requests. The trained $k$NN regressor is then queried with $\vec{\bm{s_{i}}}$ as input, after which it supplies a predictive prescription of link priorities in the form of $\vec{\bm{\beta_{i}}}$. Once the link priority vector is standardized, it is then implemented in the parameterized surrogate objective function, along with the chosen $\alpha$-profile.

\section{Numerical experiments}
\label{num_exp}

We conduct numerical experiments to assess the performance of the proposed data-driven optimization approach, i.e. the Surrogate ILP, for the online DDSP. We use the Myopic ILP as a baseline for benchmarking and we also conduct sensitivity analyses on learning parameters of the ML-enabled surrogate objective function. 
\\
\\
We conduct these experiments on the Sioux Falls network, since it is an established data set for traffic routing problems and has been used in by~\cite{levin2023branch} for the deterministic DDSP. There are 24 nodes and 76 links in this network. The length of the time horizon is $T$ = 60 minutes, with 12 intervals of 5 minutes each ($I$ = 12 and $D$ = 5 minutes). For routing and capacity tracking, the time horizon is discretized into time steps of a single minute. After 60 minutes elapse, the network still remains operational as long as the final accepted request reaches its destination. Thus the network remains active until $\max\limits_{r \in R_{\checkmark}} d_{r}$. We focus on the proof-of-concept for the proposed novel solution methodology. In line with this goal, we maintain all link capacities at 1 drone per minute, even though this is probably more restrictive than is necessary to avoid collisions. The size of the ILPs are a function of the number of variables, and assigning small capacities permits solving the problem on a large network with capacity constraints. The number of customer requests submitted in each interval follows a Poisson process with rate $\lambda$ = 100, which amounts to an expected total of 1200 requests over the course of 12 intervals. Within an interval $i$, the request submission times are uniformly distributed between the start and end of the interval $[(i-1)D, iD) $. The spatial distribution of request origins $o_{r}$ and destinations $d_{r}$ is exponential with respect to the y-coordinate, with origins dominating the bottom of the network and destinations occupying the top. For each request, the earliest permissible departure time $e_{r}$ varies uniformly between $[0,10]$ minutes of request submission time $t_{r}$. The arrival time window at the destination $[l_{r}, u_{r}]$ follows a skewed normal distribution between expected travel time and the end of the time horizon, skewed towards expected travel time. The price (profit) of each request is uniformly distributed between $[1,10]$.
\\
\\
We implemented the proposed algorithms in Python on a Windows machine with 3.19 GHz and 64 GB of RAM. We use CPLEX 22.1.0.0 to solve the test instances, and we use a solve time limit of 300 seconds for each interval of a given instance. 

\subsection{$\alpha$-profiles}

We consider three types of $\alpha$-profiles: constant, step-wise and continuous. The performance of constant profiles can be used to determine the desirable range of $\alpha_{i}$, after which temporal variations can be introduced to extend the analysis. It should be noted that even though the time intervals are discrete, a continuous $\alpha$-profile is still permissible, because $\alpha_{i}$ is only evaluated at discrete points marking the beginning of each interval. A summary of the $\alpha$-profiles experimented with is provided in Table~\ref{tab:surr_legend}. 
\\
\\
An intuitive approach consists of starting with a high value of $\alpha_{i}$ and reduce it over time. This realizes the strategy of preferentially creating link priority-based reserve capacity in the initial intervals and transitioning into profit accumulation towards the end of the time horizon. All non-constant profiles presented are non-increasing for this reason. However, constant $\alpha$-profiles have significant potential too, since the $k$NN model is adaptive and returns a different link priority vector at each interval. This implies that the priority-based slack creation is also adaptive and may relieve some of the adjustment burden from $\alpha_{i}$.

\begin{table}[t]
\centering
\renewcommand\arraystretch{1.5}
\begin{tabular}{clc}
\hline
Name     & \multicolumn{1}{c}{$\alpha$-profile composition}         & Type                        \\ \hline
SP\_CTE1 & $\alpha = 2$                                 & \multirow{6}{*}{Constant}   \\
SP\_CTE2 & $\alpha = 1.5$                               &                             \\
SP\_CTE3 & $\alpha = 1.25$                              &                             \\
SP\_CTE4 & $\alpha = 1$                                 &                             \\
SP\_CTE5 & $\alpha = 0.5$                               &                             \\
SP\_CTE6 & $\alpha = -1$                                &                             \\ \hline
SP\_STP1 & $\alpha$ = 1.5 → 1.25   | 6:6                & \multirow{6}{*}{Step-wise}  \\
SP\_STP2 & $\alpha$ = 1.5 → 1.25   | 8:4                &                             \\
SP\_STP3 & $\alpha$ = 1.5 → 1.25   → 1 | 4:4:4          &                             \\
SP\_STP4 & $\alpha$ = 1.5 → 1.25   → 1 | 7:3:2          &                             \\
SP\_STP5 & $\alpha$ = 1.25 → 1 |   6:6                  &                             \\
SP\_STP6 & $\alpha$ = 1.25 → 1 |   8:4                  &                             \\ \hline
SP\_PLY1 & $\alpha(t) = 0.01   t^{2} - 0.15545t + 1.5$  & \multirow{4}{*}{Continuous} \\
SP\_PLY2 & $\alpha(t) = -0.01   t^{2} + 0.06455t + 1.5$ &                             \\
SP\_PLY3 & $\alpha(t) = -   0.04545t + 1.5$             &                             \\
SP\_PLY4 & $\alpha(t)   =e^{-0.063t} + 0.5$             &                             \\ \hline
\end{tabular}
\caption{Summary of $\alpha$-profiles for the Surrogate ILP}
\label{tab:surr_legend}
\end{table}

\subsection{Myopic ILP policy}

We first examine the performance of the Myopic ILP policy on 5 problem instances based on the Sioux Falls network as described above. Table~\ref{tab:myopic} encapsulates the results. The Myopic ILP achieved an average profit of 3178 across all 5 instances, with the expected maximum possible profit being $p_{avg}\mathbf{E}[|\mathcal{R}|] = 6600$. The average service rate was 43.6\%. These performance results serve as a baseline for bench-marking the proposed data-driven optimization approach. The average solve-time for an interval across all 5 instances was 123.3 seconds. We observed that if an optimal solution was reached, then the solve-time was below 50 seconds, but for some intervals the solver reached the time limit of 300 seconds and gave the best feasible solution it could find. The latter case occurred in the initial three to four intervals, presumably because the network traffic was still evolving and had not reached complete saturation. In such a situation, many idle requests have not yet been dispatched and incoming requests continue to arrive at each interval, making it difficult to reach optimality. As idle requests begin to enter the network in greater numbers and network traffic increases, the pending request load subsides and optimality can be reached more easily. 

\begin{table}[t]
\centering
\renewcommand\arraystretch{1.5}
\begin{tabular}{cccc}
\hline
Instance & Total Requests & Service Rate & Profit \\ \hline
0        & 1203           & 41.6\%       & 3076   \\
1        & 1125           & 45.5\%       & 3249   \\
2        & 1260           & 41.7\%       & 3260   \\
3        & 1184           & 44.5\%       & 3216   \\
4        & 1147           & 44.0\%       & 3089   \\ \hline
\end{tabular}
\caption{Performance of the Myopic ILP}
  \label{tab:myopic}
\end{table}

\subsection{Training $k$NN multi-output regressor}

\begin{figure}[H]
\centering
\hspace{-2cm}
\begin{subfigure}{0.46\textwidth}
    \includegraphics[width=\textwidth]{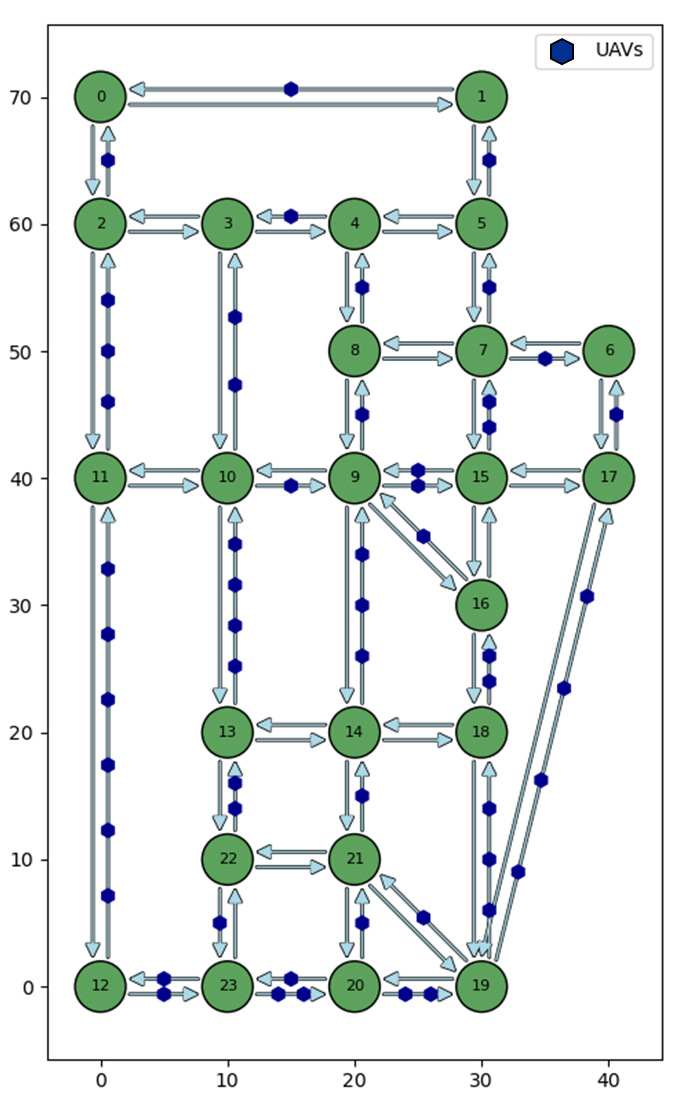}
    \caption{Input: Link Occupancy vector $\bm{\vec{s_{i}}} = \begin{bmatrix}
    s_{i, l}
\end{bmatrix}^{l \in A}$}
    \label{fig:first}
\end{subfigure}
\hspace{0cm}
\begin{subfigure}{0.46\textwidth}
    \includegraphics[width=1.25\textwidth]{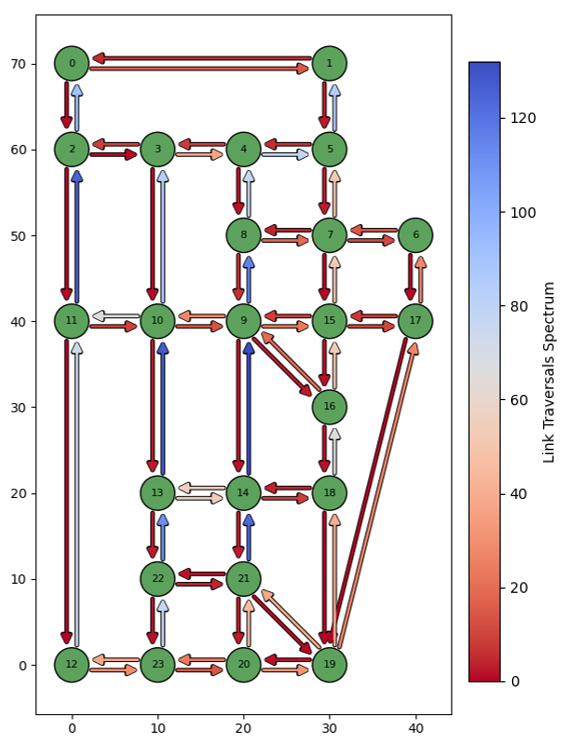}
    \caption{Target Output: Link Priority vector $\bm{\vec{\beta_{i}}} = \begin{bmatrix}
   \beta_{i, l}
\end{bmatrix}^{l \in A}$}
    \label{fig:second}
\end{subfigure}  

\caption{(a) Input and (b) Target output vectors for training interval 500 in the training scheme. \textit{Note: drones occupying the same link in the network need not be equidistant  but have been shown as such for simplicity.}}
\label{fig:visual_dem}
\end{figure}

The training instance contains 2000 intervals, instead of the usual 12 for test instances. The lookahead simulation contained 5 virtual intervals. We ran the training scheme described in Algorithm~\ref{alg:beta-train} on this instance and obtained the training dataset where the input feature set is $\bm{\vec{S}}$ and the output target set is $\bm{\vec{B}}$. We employed the reservation heuristic to handle the size of the training instance. The time-expanded dual graph $G^{T}$ is constructed up to a time layer $t_{\max} + \Delta$ where $t_{\max} = \max\limits_{r \in \mathcal{R}} d_{r}$ of the training requests. The purpose of the buffer $\Delta$ is to account for the virtual requests being created in the final interval, some of whose arrival windows may exceed $t_{\max}$.
\\
\\
The completed training dataset is passed through a $k$NN multi-output regressor. We implemented $k$NN through the Scikit-learn package in Python. As mentioned previously, the collection vectors $\bm{\vec{S}}$ and $\bm{\vec{B}}$ contain the link occupancy vectors $\bm{\vec{s_{i}}} = \begin{bmatrix}
    s_{i, l}
\end{bmatrix}^{l \in A}$ and link priority vectors $\bm{\vec{\beta_{i}}} = \begin{bmatrix}
   \beta_{i, l}
\end{bmatrix}^{l \in A}$ respectively for all training intervals. During the execution of the training scheme, we converted the link occupancy vector $\bm{\vec{s_{i}}}$ into a 2-D node adjacency matrix $\bm{N_{s}}$ where the row and column indices of the matrix are the nodes of the network. Each entry $s_{m,n}$ of this matrix represents a directed link from node $m$ to node $n$ along link $\omega(m,n)$, and the value of the entry equals the link occupancy. This transformation is useful because it encodes the state of the network in a more information-rich format. Apart from storing link occupancy, the node adjacency matrix also conveys the network topology (for example, two separate non-zero entries in the same row $m$ represent two outgoing links with $m$ as the common source). With $\bm{N_{s}}$ as input, the $k$NN regressor is trained to estimate the link priorities as output. The number of neighbours considered for each input is set to 60, and the distance metric is set to Euclidean distance.
\\
\\
Figure~\ref{fig:first} demonstrates one pair of input and target output vectors in the training scheme (training interval 500). The input is the snapshot of the network at the beginning of training interval 500, depicting the number of drones (UAVs) on each link. After simulating 5 virtual intervals, the link traversals of the accepted virtual requests give rise to the desired target output as shown in Figure~\ref{fig:second}. This $k$NN model is trained on such input and target output pairs for the entire duration of the training scheme. We note that the spatial distribution of request origins and destinations influences the link priorities. Since origins dominate the bottom portion of the network and destinations dominate the top, links directed southwards (top to bottom) are of least priority. Links within the central portion of the network, which carry traffic in the opposite direction (bottom to top) are of high priority. Some links responsible for lateral movement are also of moderate to high importance. 
    
\subsection{Constant $\alpha_{i}$}

We conduct a sensitivity analysis to quantify the impact of the $\alpha$-profile onto the performance of the Surrogate ILP policy. Since $\alpha_i = 0 \ \forall i \in \{1, \ldots, I\}$ is equivalent to the Myopic ILP, we perturbed $\alpha_i$ above and below 0 and observed the effect on cumulative profit acquired. The cumulative profit of the Surrogate ILP is used to determine the Profit Gap \% with respect to the Myopic ILP and serves as the main performance metric. We similarly recorded the cumulative service (total requests accepted) and Service Gap \%. The performance of the Surrogate ILP policy with constant $\alpha$-profiles are displayed in Table~\ref{tab:ctalpha}. 
\\
\\
As expected, a fully negative $\alpha$-profile (SP\_CTE6) performed worse than the Myopic ILP, because the surrogate objective seeks to minimize slack and force traffic though the network in a counter-productive manner. The influence of the link priorities $\beta_{i,l}$ further worsened the solution, because high priority links are overloaded more in proportion to their priority. In contrast, a positive alpha profile performed better than the benchmark, but the profit gap decreases if $\alpha$ is exceedingly large. This can be explained by observing that too large an alpha maintained over the entire time horizon drains the network excessively and prevents profitable requests from being routed successfully. For the given network and the demand scenario presented, constant $\alpha$ profiles of 1.5, 1.25 and 1 performed the best, with an average gap across the 5 problem instances of 9.43\%, 10.62\% and 9.73\% respectively. For these three $\alpha$ profiles, the average solve-times for each interval across the 5 problem instances are 148.5, 150.0, 142.4 seconds respectively. A similar pattern to the myopic ILP emerged, whereby initial three to four intervals deplete the solve-time limit but thereafter optimality is reached well before it.
\\
\\
We observed that the Surrogate ILP policies served slightly fewer requests than the Myopic ILP policy, but attained a larger profit. The Myopic ILP behaved greedily -- optimal in the current interval (for initial intervals) but negligent of the future -- and lost out on potentially profitable requests in later intervals. Instead, the Surrogate ILPs maintained reserve capacity on crucial links which enabled the acceptance of high profit requests no matter in which interval they arrived. It is possible however, to create excessive reserves that begin to hinder traffic through the network, as evidenced by the inferior performance (relative to other Surrogate ILPs) of SP\_CTE1. With $\alpha_{i} = 2$, SP\_CTE1 had the lowest service rate among all Surrogate ILPs but was also the worst performing policy that still outperformed the Myopic ILP policy.

\begin{table}[t]
\centering
\resizebox{\textwidth}{!}{
\renewcommand\arraystretch{1.5}
\begin{tabular}{ccccccccccc}
\hline
\multirow{2}{*}{Instance} & \multirow{2}{*}{Total Arrived} & \multicolumn{2}{c}{Myopic ILP}                & \multirow{2}{*}{$\alpha$-profile} & \multirow{2}{*}{Service Rate} & \multirow{2}{*}{Profit} & \multirow{2}{*}{Profit Gap} & \multirow{2}{*}{Profit Gap \%} & \multirow{2}{*}{Service Gap} & \multirow{2}{*}{Service Gap \%} \\ \cline{3-4}
                          &                                & \multicolumn{1}{l}{Service Rate} & Profit                &                            &                               &                         &                             &                                &                              &                                 \\ \hline
\multirow{6}{*}{0}        & \multirow{6}{*}{1203}          & \multirow{6}{*}{41.6\%}          & \multirow{6}{*}{3076} & SP\_CTE1                   & 37.0\%                        & 3255                    & 179                         & 5.8\%                          & -56                          & -4.7\%                          \\
                          &                                &                                  &                       & SP\_CTE2                   & 40.6\%                        & 3497                    & 421                         & 13.7\%                         & -13                          & -1.0\%                          \\
                          &                                &                                  &                       & SP\_CTE3                   & 40.9\%                        & 3474                    & 398                         & 12.9\%                         & -9                           & -0.7\%                          \\
                          &                                &                                  &                       & SP\_CTE4                   & 41.6\%                        & 3413                    & 337                         & 11.0\%                         & -1                           & 0.0\%                           \\
                          &                                &                                  &                       & SP\_CTE5                   & 42.6\%                        & 3308                    & 232                         & 7.5\%                          & 11                           & 1.0\%                           \\
                          &                                &                                  &                       & SP\_CTE6                   & 41.5\%                        & 3062                    & -14                         & -0.5\%                         & -2                           & -0.1\%                          \\ \hline
\multirow{6}{*}{1}        & \multirow{6}{*}{1125}          & \multirow{6}{*}{45.5\%}          & \multirow{6}{*}{3249} & SP\_CTE1                   & 40.1\%                        & 3290                    & 41                          & 1.3\%                          & -61                          & -5.4\%                          \\
                          &                                &                                  &                       & SP\_CTE2                   & 43.8\%                        & 3534                    & 285                         & 8.8\%                          & -19                          & -1.7\%                          \\
                          &                                &                                  &                       & SP\_CTE3                   & 45.0\%                        & 3590                    & 341                         & 10.5\%                         & -6                           & -0.5\%                          \\
                          &                                &                                  &                       & SP\_CTE4                   & 45.6\%                        & 3542                    & 293                         & 9.0\%                          & 1                            & 0.1\%                           \\
                          &                                &                                  &                       & SP\_CTE5                   & 46.4\%                        & 3409                    & 160                         & 4.9\%                          & 10                           & 0.9\%                           \\
                          &                                &                                  &                       & SP\_CTE6                   & 44.6\%                        & 3141                    & -108                        & -3.3\%                         & -10                          & -0.9\%                          \\ \hline
\multirow{6}{*}{2}        & \multirow{6}{*}{1260}          & \multirow{6}{*}{41.7\%}          & \multirow{6}{*}{3260} & SP\_CTE1                   & 36.4\%                        & 3328                    & 68                          & 2.1\%                          & -67                          & -5.3\%                          \\
                          &                                &                                  &                       & SP\_CTE2                   & 39.7\%                        & 3563                    & 303                         & 9.3\%                          & -26                          & -2.0\%                          \\
                          &                                &                                  &                       & SP\_CTE3                   & 40.6\%                        & 3577                    & 317                         & 9.7\%                          & -14                          & -1.1\%                          \\
                          &                                &                                  &                       & SP\_CTE4                   & 41.3\%                        & 3564                    & 304                         & 9.3\%                          & -6                           & -0.4\%                          \\
                          &                                &                                  &                       & SP\_CTE5                   & 41.9\%                        & 3401                    & 141                         & 4.3\%                          & 2                            & 0.2\%                           \\
                          &                                &                                  &                       & SP\_CTE6                   & 41.4\%                        & 3222                    & -38                         & -1.2\%                         & -4                           & -0.3\%                          \\ \hline
\multirow{6}{*}{3}        & \multirow{6}{*}{1184}          & \multirow{6}{*}{44.5\%}          & \multirow{6}{*}{3216} & SP\_CTE1                   & 37.8\%                        & 3280                    & 64                          & 2.0\%                          & -80                          & -6.8\%                          \\
                          &                                &                                  &                       & SP\_CTE2                   & 40.7\%                        & 3443                    & 227                         & 7.1\%                          & -45                          & -3.8\%                          \\
                          &                                &                                  &                       & SP\_CTE3                   & 42.1\%                        & 3502                    & 286                         & 8.9\%                          & -29                          & -2.4\%                          \\
                          &                                &                                  &                       & SP\_CTE4                   & 43.3\%                        & 3515                    & 299                         & 9.3\%                          & -14                          & -1.2\%                          \\
                          &                                &                                  &                       & SP\_CTE5                   & 43.6\%                        & 3278                    & 62                          & 1.9\%                          & -11                          & -0.9\%                          \\
                          &                                &                                  &                       & SP\_CTE6                   & 43.7\%                        & 3195                    & -21                         & -0.7\%                         & -10                          & -0.8\%                          \\ \hline
\multirow{6}{*}{4}        & \multirow{6}{*}{1147}          & \multirow{6}{*}{44.0\%}          & \multirow{6}{*}{3089} & SP\_CTE1                   & 37.1\%                        & 3148                    & 59                          & 1.9\%                          & -79                          & -6.9\%                          \\
                          &                                &                                  &                       & SP\_CTE2                   & 41.0\%                        & 3346                    & 257                         & 8.3\%                          & -35                          & -3.0\%                          \\
                          &                                &                                  &                       & SP\_CTE3                   & 43.0\%                        & 3430                    & 341                         & 11.0\%                         & -12                          & -1.0\%                          \\
                          &                                &                                  &                       & SP\_CTE4                   & 43.0\%                        & 3400                    & 311                         & 10.1\%                         & -12                          & -1.0\%                          \\
                          &                                &                                  &                       & SP\_CTE5                   & 43.9\%                        & 3590                    & 501                         & 16.2\%                         & -1                           & -0.1\%                          \\
                          &                                &                                  &                       & SP\_CTE6                   & 44.1\%                        & 2980                    & -109                        & -3.5\%                         & 1                            & 0.1\%                           \\ \hline
\end{tabular}
}
\caption{Surrogate ILP policy with constant $\alpha$-profile vs. Myopic ILP policy}
\label{tab:ctalpha}
\end{table}

\subsection{Step-wise \& continuous $\alpha_{i}$}

Having demonstrated the persistent superiority of SP\_CTE2, SP\_CTE3 and SP\_CTE4, we explore step-wise changes in the $\alpha$-profile with 1.5, 1.25 and 1 as pivots. We named these Surrogate ILPs SP\_STP1 through SP\_STP6. For example, SP\_STP1 has the $\alpha$-profile $1.5 \rightarrow 1.25 \ | \ 6:6$, which depicts 6 intervals of $\alpha_{i} = 1.5$ followed by 6 intervals of $\alpha_{i} = 1.25$. Similarly, SP\_STP2 has the $\alpha$-profile $1.5 \rightarrow 1.25 \ | \ 8:4$, which depicts 8 intervals of $\alpha_{i} = 1.5$ followed by 4 intervals of $\alpha_{i} = 1.25$. The subsequent two step-wise profiles move from 1.5 all the way to 1, and the final two go from 1.25 to 1. In these 3 pairs of Surrogate ILP policies, the latter one captures a delayed transition to a lower $\alpha_{i}$.
\\
\\
Table~\ref{tab:swalpha} records the performance of this class of $\alpha$-profiles. We observe a higher quality solution from the group of step-wise $\alpha$-profiles than from the constant profiles. However, there is no clear superiority of one step-wise profile over another. We find that only SP\_STP5 emerged as the best in 2 instances, while the remaining 3 instances are each dominated by adifferent step profile. As with SP\_CTE2, SP\_CTE3 and SP\_CTE4, the service rate is less than that of the Myopic ILP policy but the cumulative profit is greater.
\\
\\
The link priority vector returned by the $k$NN multi-output regressor is adaptive and different in each interval, since it is extremely unlikely for the input snapshot vector of link occupancy to be identical in two consecutive intervals. For this reason, it is difficult to establish a prevailing improvement over constant $\alpha$-profiles with step-wise or continuous profiles. By continuous $\alpha_{i}$ we imply a temporal function for the $\alpha$-profile which is evaluated at discrete time-steps (the start of each interval). We approximated the step-wise algorithms SP\_STP3 and SP\_STP4 as continuous functions of time. These algorithms (SP\_PLY1 to SP\_PLY4) outperformed the Myopic ILP as well, as shown in Table~\ref{tab:plalpha}.
\\
\\
In Table~\ref{tab:cte_vs_cnt}, we present a comparison of the step-wise and continuous $\alpha$-profiles with the 3 best constant profiles. No sustained performance gap is observed, indicating that it is difficult to generate consistent improvement upon the adaptive $k$NN component, provided the value of $\alpha$ is not arbitrarily chosen. In Table~\ref{tab:htmap} we summarize the profits and profit gaps (relative to the Myopic ILP) of all surrogate objective functions in the form of a heat map.

\begin{table}[t]
\centering
\resizebox{\textwidth}{!}{
\renewcommand\arraystretch{1.5}
\begin{tabular}{ccccccccccc}
\hline
\multirow{2}{*}{Instance} & \multirow{2}{*}{Total Arrived} & \multicolumn{2}{c}{Myopic ILP}                & \multirow{2}{*}{$\alpha$-profile} & \multirow{2}{*}{Service Rate} & \multirow{2}{*}{Profit} & \multirow{2}{*}{Profit Gap} & \multirow{2}{*}{Profit Gap \%} & \multirow{2}{*}{Service Gap} & \multirow{2}{*}{Service Gap \%} \\ \cline{3-4}
                          &                                & \multicolumn{1}{l}{Service Rate} & Profit                &                            &                               &                         &                             &                                &                              &                                 \\ \hline
\multirow{6}{*}{0}        & \multirow{6}{*}{1203}          & \multirow{6}{*}{41.6\%}          & \multirow{6}{*}{3076} & SP\_STP1                   & 40.2\%                        & 3444                    & 368                         & 12.0\%                         & -17                          & -1.4\%                          \\
                          &                                &                                  &                       & SP\_STP2                   & 40.5\%                        & 3498                    & 422                         & 13.7\%                         & -14                          & -1.1\%                          \\
                          &                                &                                  &                       & SP\_STP3                   & 41.3\%                        & 3501                    & 425                         & 13.8\%                         & -4                           & -0.3\%                          \\
                          &                                &                                  &                       & SP\_STP4                   & 40.2\%                        & 3481                    & 405                         & 13.2\%                         & -17                          & -1.4\%                          \\
                          &                                &                                  &                       & SP\_STP5                   & 41.6\%                        & 3511                    & 435                         & 14.1\%                         & -1                           & 0.0\%                           \\
                          &                                &                                  &                       & SP\_STP6                   & 41.1\%                        & 3449                    & 373                         & 12.1\%                         & -7                           & -0.5\%                          \\ \hline
\multirow{6}{*}{1}        & \multirow{6}{*}{1125}          & \multirow{6}{*}{45.5\%}          & \multirow{6}{*}{3249} & SP\_STP1                   & 43.8\%                        & 3568                    & 319                         & 9.8\%                          & -19                          & -1.7\%                          \\
                          &                                &                                  &                       & SP\_STP2                   & 44.4\%                        & 3584                    & 335                         & 10.3\%                         & -12                          & -1.1\%                          \\
                          &                                &                                  &                       & SP\_STP3                   & 43.6\%                        & 3485                    & 236                         & 7.3\%                          & -21                          & -1.9\%                          \\
                          &                                &                                  &                       & SP\_STP4                   & 44.7\%                        & 3588                    & 339                         & 10.4\%                         & -9                           & -0.8\%                          \\
                          &                                &                                  &                       & SP\_STP5                   & 45.2\%                        & 3547                    & 298                         & 9.2\%                          & -3                           & -0.3\%                          \\
                          &                                &                                  &                       & SP\_STP6                   & 44.4\%                        & 3517                    & 268                         & 8.2\%                          & -12                          & -1.1\%                          \\ \hline
\multirow{6}{*}{2}        & \multirow{6}{*}{1260}          & \multirow{6}{*}{41.7\%}          & \multirow{6}{*}{3260} & SP\_STP1                   & 40.5\%                        & 3591                    & 331                         & 10.2\%                         & -16                          & -1.2\%                          \\
                          &                                &                                  &                       & SP\_STP2                   & 40.4\%                        & 3614                    & 354                         & 10.9\%                         & -17                          & -1.3\%                          \\
                          &                                &                                  &                       & SP\_STP3                   & 40.8\%                        & 3570                    & 310                         & 9.5\%                          & -12                          & -0.9\%                          \\
                          &                                &                                  &                       & SP\_STP4                   & 40.2\%                        & 3568                    & 308                         & 9.4\%                          & -19                          & -1.5\%                          \\
                          &                                &                                  &                       & SP\_STP5                   & 41.1\%                        & 3582                    & 322                         & 9.9\%                          & -8                           & -0.6\%                          \\
                          &                                &                                  &                       & SP\_STP6                   & 41.5\%                        & 3625                    & 365                         & 11.2\%                         & -3                           & -0.2\%                          \\ \hline
\multirow{6}{*}{3}        & \multirow{6}{*}{1184}          & \multirow{6}{*}{44.5\%}          & \multirow{6}{*}{3216} & SP\_STP1                   & 41.9\%                        & 3536                    & 320                         & 10.0\%                         & -31                          & -2.6\%                          \\
                          &                                &                                  &                       & SP\_STP2                   & 41.6\%                        & 3492                    & 276                         & 8.6\%                          & -34                          & -2.9\%                          \\
                          &                                &                                  &                       & SP\_STP3                   & 42.1\%                        & 3470                    & 254                         & 7.9\%                          & -29                          & -2.4\%                          \\
                          &                                &                                  &                       & SP\_STP4                   & 42.4\%                        & 3552                    & 336                         & 10.4\%                         & -25                          & -2.1\%                          \\
                          &                                &                                  &                       & SP\_STP5                   & 43.6\%                        & 3556                    & 340                         & 10.6\%                         & -11                          & -0.9\%                          \\
                          &                                &                                  &                       & SP\_STP6                   & 42.5\%                        & 3470                    & 254                         & 7.9\%                          & -24                          & -2.0\%                          \\ \hline
\multirow{6}{*}{4}        & \multirow{6}{*}{1147}          & \multirow{6}{*}{44.0\%}          & \multirow{6}{*}{3089} & SP\_STP1                   & 42.2\%                        & 3431                    & 342                         & 11.1\%                         & -21                          & -1.8\%                          \\
                          &                                &                                  &                       & SP\_STP2                   & 41.7\%                        & 3395                    & 306                         & 9.9\%                          & -27                          & -2.3\%                          \\
                          &                                &                                  &                       & SP\_STP3                   & 42.9\%                        & 3460                    & 371                         & 12.0\%                         & -13                          & -1.1\%                          \\
                          &                                &                                  &                       & SP\_STP4                   & 41.8\%                        & 3382                    & 293                         & 9.5\%                          & -26                          & -2.2\%                          \\
                          &                                &                                  &                       & SP\_STP5                   & 43.0\%                        & 3412                    & 323                         & 10.5\%                         & -12                          & -1.0\%                          \\
                          &                                &                                  &                       & SP\_STP6                   & 42.6\%                        & 3418                    & 329                         & 10.7\%                         & -16                          & -1.4\%                          \\ \hline
\end{tabular}
}
\caption{Surrogate ILP policy with step-wise $\alpha$-profile vs. Myopic ILP policy}
\label{tab:swalpha}
\end{table}

\begin{table}[t]
\centering
\resizebox{\textwidth}{!}{
\renewcommand\arraystretch{1.5}
\begin{tabular}{ccccccccccc}
\hline
\multirow{2}{*}{Instance} & \multirow{2}{*}{Total Arrived} & \multicolumn{2}{c}{Myopic ILP}                & \multirow{2}{*}{$\alpha$-profile} & \multirow{2}{*}{Service Rate} & \multirow{2}{*}{Profit} & \multirow{2}{*}{Profit Gap} & \multirow{2}{*}{Profit Gap \%} & \multirow{2}{*}{Service Gap} & \multirow{2}{*}{Service Gap \%} \\ \cline{3-4}
                          &                                & \multicolumn{1}{l}{Service Rate} & Profit                &                            &                               &                         &                             &                                &                              &                                 \\ \hline
\multirow{4}{*}{0}        & \multirow{4}{*}{1203}          & \multirow{4}{*}{41.6\%}          & \multirow{4}{*}{3076} & SP\_PLY1                   & 41.1\%                        & 3461                    & 385                         & 12.5\%                         & -7                           & -0.5\%                          \\
                          &                                &                                  &                       & SP\_PLY2                   & 40.6\%                        & 3502                    & 426                         & 13.8\%                         & -13                          & -1.0\%                          \\
                          &                                &                                  &                       & SP\_PLY3                   & 41.7\%                        & 3542                    & 466                         & 15.1\%                         & 1                            & 0.1\%                           \\
                          &                                &                                  &                       & SP\_PLY4                   & 40.9\%                        & 3458                    & 382                         & 12.4\%                         & -9                           & -0.7\%                          \\ \hline
\multirow{4}{*}{1}        & \multirow{4}{*}{1125}          & \multirow{4}{*}{45.5\%}          & \multirow{4}{*}{3249} & SP\_PLY1                   & 44.6\%                        & 3514                    & 265                         & 8.2\%                          & -10                          & -0.9\%                          \\
                          &                                &                                  &                       & SP\_PLY2                   & 44.3\%                        & 3571                    & 322                         & 9.9\%                          & -14                          & -1.2\%                          \\
                          &                                &                                  &                       & SP\_PLY3                   & 45.0\%                        & 3562                    & 313                         & 9.6\%                          & -6                           & -0.5\%                          \\
                          &                                &                                  &                       & SP\_PLY4                   & 44.7\%                        & 3530                    & 281                         & 8.6\%                          & -9                           & -0.8\%                          \\ \hline
\multirow{4}{*}{2}        & \multirow{4}{*}{1260}          & \multirow{4}{*}{41.7\%}          & \multirow{4}{*}{3260} & SP\_PLY1                   & 40.9\%                        & 3540                    & 280                         & 8.6\%                          & -11                          & -0.8\%                          \\
                          &                                &                                  &                       & SP\_PLY2                   & 40.9\%                        & 3648                    & 388                         & 11.9\%                         & -11                          & -0.8\%                          \\
                          &                                &                                  &                       & SP\_PLY3                   & 40.6\%                        & 3572                    & 312                         & 9.6\%                          & -14                          & -1.1\%                          \\
                          &                                &                                  &                       & SP\_PLY4                   & 40.8\%                        & 3576                    & 316                         & 9.7\%                          & -12                          & -0.9\%                          \\ \hline
\multirow{4}{*}{3}        & \multirow{4}{*}{1184}          & \multirow{4}{*}{44.5\%}          & \multirow{4}{*}{3216} & SP\_PLY1                   & 43.2\%                        & 3501                    & 285                         & 8.9\%                          & -16                          & -1.3\%                          \\
                          &                                &                                  &                       & SP\_PLY2                   & 42.2\%                        & 3544                    & 328                         & 10.2\%                         & -27                          & -2.3\%                          \\
                          &                                &                                  &                       & SP\_PLY3                   & 42.8\%                        & 3546                    & 330                         & 10.3\%                         & -20                          & -1.7\%                          \\
                          &                                &                                  &                       & SP\_PLY4                   & 42.8\%                        & 3550                    & 334                         & 10.4\%                         & -20                          & -1.7\%                          \\ \hline
\multirow{4}{*}{4}        & \multirow{4}{*}{1147}          & \multirow{4}{*}{44.0\%}          & \multirow{4}{*}{3089} & SP\_PLY1                   & 43.8\%                        & 3436                    & 347                         & 11.2\%                         & -3                           & -0.2\%                          \\
                          &                                &                                  &                       & SP\_PLY2                   & 41.8\%                        & 3379                    & 290                         & 9.4\%                          & -26                          & -2.2\%                          \\
                          &                                &                                  &                       & SP\_PLY3                   & 42.5\%                        & 3411                    & 322                         & 10.4\%                         & -18                          & -1.5\%                          \\
                          &                                &                                  &                       & SP\_PLY4                   & 42.9\%                        & 3453                    & 364                         & 11.8\%                         & -13                          & -1.1\%                          \\ \hline
\end{tabular}
}
\caption{Surrogate ILP policy with continuous $\alpha$-profile vs. Myopic ILP policy}
  \label{tab:plalpha}
\end{table}

\begin{table}[t]
\centering
\resizebox{\textwidth}{!}{
\renewcommand\arraystretch{1.5}
\begin{tabular}{clccccccccccc}
\hline
\multicolumn{2}{c}{}                                            &                            & \multicolumn{10}{c}{Step-wise and Continuous $\alpha$-profiles}                                                                                                                                                                                                                                                                                \\ \cline{4-13} 
\multicolumn{2}{c}{\multirow{-2}{*}{Constant $\alpha$-profile}} & \multirow{-2}{*}{Instance} & SP\_STP1                       & SP\_STP2                       & SP\_STP3                       & SP\_STP4                       & SP\_STP5                       & SP\_STP6                       & SP\_PLY1                       & SP\_PLY2                       & SP\_PLY3                       & SP\_PLY4                       \\ \hline
\multicolumn{2}{c}{}                                            & 0                          & \cellcolor[HTML]{FA9F75}-1.5\% & \cellcolor[HTML]{FEDB81}0.0\%  & \cellcolor[HTML]{FEDE81}0.1\%  & \cellcolor[HTML]{FDC87D}-0.5\% & \cellcolor[HTML]{FEE983}0.4\%  & \cellcolor[HTML]{FBA576}-1.4\% & \cellcolor[HTML]{FBB279}-1.0\% & \cellcolor[HTML]{FEDF81}0.1\%  & \cellcolor[HTML]{D8E082}1.3\%  & \cellcolor[HTML]{FBAF78}-1.1\% \\
\multicolumn{2}{c}{}                                            & 1                          & \cellcolor[HTML]{E7E483}1.0\%  & \cellcolor[HTML]{D2DE82}1.4\%  & \cellcolor[HTML]{FBA476}-1.4\% & \cellcolor[HTML]{CCDD82}1.5\%  & \cellcolor[HTML]{FEE883}0.4\%  & \cellcolor[HTML]{FDC77D}-0.5\% & \cellcolor[HTML]{FCC47C}-0.6\% & \cellcolor[HTML]{E3E383}1.0\%  & \cellcolor[HTML]{EFE784}0.8\%  & \cellcolor[HTML]{FDD67F}-0.1\% \\
\multicolumn{2}{c}{}                                            & 2                          & \cellcolor[HTML]{EFE784}0.8\%  & \cellcolor[HTML]{D1DE82}1.4\%  & \cellcolor[HTML]{FEE282}0.2\%  & \cellcolor[HTML]{FEDF81}0.1\%  & \cellcolor[HTML]{FBEA84}0.5\%  & \cellcolor[HTML]{C2DA81}1.7\%  & \cellcolor[HTML]{FCC17B}-0.6\% & \cellcolor[HTML]{A4D17F}2.4\%  & \cellcolor[HTML]{FEE482}0.3\%  & \cellcolor[HTML]{FEE883}0.4\%  \\
\multicolumn{2}{c}{}                                            & 3                          & \cellcolor[HTML]{96CD7E}2.7\%  & \cellcolor[HTML]{D1DE82}1.4\%  & \cellcolor[HTML]{EFE784}0.8\%  & \cellcolor[HTML]{80C77D}3.2\%  & \cellcolor[HTML]{7BC57D}3.3\%  & \cellcolor[HTML]{EFE784}0.8\%  & \cellcolor[HTML]{C5DB81}1.7\%  & \cellcolor[HTML]{8BCA7E}2.9\%  & \cellcolor[HTML]{88C97E}3.0\%  & \cellcolor[HTML]{83C87D}3.1\%  \\
\multicolumn{2}{c}{\multirow{-5}{*}{SP\_CTE2}}                  & 4                          & \cellcolor[HTML]{9DCF7F}2.5\%  & \cellcolor[HTML]{CFDE82}1.5\%  & \cellcolor[HTML]{75C47D}3.4\%  & \cellcolor[HTML]{E1E383}1.1\%  & \cellcolor[HTML]{B8D780}2.0\%  & \cellcolor[HTML]{AFD480}2.2\%  & \cellcolor[HTML]{96CD7E}2.7\%  & \cellcolor[HTML]{E5E483}1.0\%  & \cellcolor[HTML]{B9D780}1.9\%  & \cellcolor[HTML]{7FC67D}3.2\%  \\ \hline
\multicolumn{2}{c}{}                                            & 0                          & \cellcolor[HTML]{FCB87A}-0.9\% & \cellcolor[HTML]{F3E884}0.7\%  & \cellcolor[HTML]{EFE784}0.8\%  & \cellcolor[HTML]{FEE282}0.2\%  & \cellcolor[HTML]{E2E383}1.1\%  & \cellcolor[HTML]{FCBE7B}-0.7\% & \cellcolor[HTML]{FDCB7E}-0.4\% & \cellcolor[HTML]{EEE683}0.8\%  & \cellcolor[HTML]{B8D780}2.0\%  & \cellcolor[HTML]{FDC87D}-0.5\% \\
\multicolumn{2}{c}{}                                            & 1                          & \cellcolor[HTML]{FCC27C}-0.6\% & \cellcolor[HTML]{FDD37F}-0.2\% & \cellcolor[HTML]{F8696B}-2.9\% & \cellcolor[HTML]{FDD880}-0.1\% & \cellcolor[HTML]{FBAB77}-1.2\% & \cellcolor[HTML]{F98B71}-2.0\% & \cellcolor[HTML]{F98871}-2.1\% & \cellcolor[HTML]{FDC57C}-0.5\% & \cellcolor[HTML]{FCBC7A}-0.8\% & \cellcolor[HTML]{FA9974}-1.7\% \\
\multicolumn{2}{c}{}                                            & 2                          & \cellcolor[HTML]{FEE983}0.4\%  & \cellcolor[HTML]{E3E383}1.0\%  & \cellcolor[HTML]{FDD27F}-0.2\% & \cellcolor[HTML]{FDD07E}-0.3\% & \cellcolor[HTML]{FEDF81}0.1\%  & \cellcolor[HTML]{D5DF82}1.3\%  & \cellcolor[HTML]{FBB279}-1.0\% & \cellcolor[HTML]{B7D780}2.0\%  & \cellcolor[HTML]{FDD57F}-0.1\% & \cellcolor[HTML]{FED980}0.0\%  \\
\multicolumn{2}{c}{}                                            & 3                          & \cellcolor[HTML]{E6E483}1.0\%  & \cellcolor[HTML]{FDCF7E}-0.3\% & \cellcolor[HTML]{FCB679}-0.9\% & \cellcolor[HTML]{D1DE82}1.4\%  & \cellcolor[HTML]{CCDD82}1.5\%  & \cellcolor[HTML]{FCB679}-0.9\% & \cellcolor[HTML]{FED980}0.0\%  & \cellcolor[HTML]{DCE182}1.2\%  & \cellcolor[HTML]{D9E082}1.3\%  & \cellcolor[HTML]{D4DF82}1.4\%  \\
\multicolumn{2}{c}{\multirow{-5}{*}{SP\_CTE3}}                  & 4                          & \cellcolor[HTML]{FEDB81}0.0\%  & \cellcolor[HTML]{FBB279}-1.0\% & \cellcolor[HTML]{EBE583}0.9\%  & \cellcolor[HTML]{FBA476}-1.4\% & \cellcolor[HTML]{FDC67C}-0.5\% & \cellcolor[HTML]{FDCC7E}-0.3\% & \cellcolor[HTML]{FEE182}0.2\%  & \cellcolor[HTML]{FBA075}-1.5\% & \cellcolor[HTML]{FCC47C}-0.6\% & \cellcolor[HTML]{F4E884}0.7\%  \\ \hline
\multicolumn{2}{c}{}                                            & 0                          & \cellcolor[HTML]{E9E583}0.9\%  & \cellcolor[HTML]{9FD07F}2.5\%  & \cellcolor[HTML]{9BCF7F}2.6\%  & \cellcolor[HTML]{B7D680}2.0\%  & \cellcolor[HTML]{8ECB7E}2.9\%  & \cellcolor[HTML]{E2E383}1.1\%  & \cellcolor[HTML]{D2DE82}1.4\%  & \cellcolor[HTML]{9ACE7F}2.6\%  & \cellcolor[HTML]{63BE7B}3.8\%  & \cellcolor[HTML]{D6E082}1.3\%  \\
\multicolumn{2}{c}{}                                            & 1                          & \cellcolor[HTML]{F1E784}0.7\%  & \cellcolor[HTML]{DCE182}1.2\%  & \cellcolor[HTML]{FA9C74}-1.6\% & \cellcolor[HTML]{D7E082}1.3\%  & \cellcolor[HTML]{FEDF81}0.1\%  & \cellcolor[HTML]{FCBF7B}-0.7\% & \cellcolor[HTML]{FCBB7A}-0.8\% & \cellcolor[HTML]{EDE683}0.8\%  & \cellcolor[HTML]{F9EA84}0.6\%  & \cellcolor[HTML]{FDCD7E}-0.3\% \\
\multicolumn{2}{c}{}                                            & 2                          & \cellcolor[HTML]{F0E784}0.8\%  & \cellcolor[HTML]{D2DE82}1.4\%  & \cellcolor[HTML]{FEE082}0.2\%  & \cellcolor[HTML]{FEDE81}0.1\%  & \cellcolor[HTML]{FCEA84}0.5\%  & \cellcolor[HTML]{C4DA81}1.7\%  & \cellcolor[HTML]{FCC07B}-0.7\% & \cellcolor[HTML]{A6D27F}2.4\%  & \cellcolor[HTML]{FEE382}0.2\%  & \cellcolor[HTML]{FEE783}0.3\%  \\
\multicolumn{2}{c}{}                                            & 3                          & \cellcolor[HTML]{F8E984}0.6\%  & \cellcolor[HTML]{FCC17B}-0.7\% & \cellcolor[HTML]{FBA877}-1.3\% & \cellcolor[HTML]{E2E383}1.1\%  & \cellcolor[HTML]{DDE283}1.2\%  & \cellcolor[HTML]{FBA877}-1.3\% & \cellcolor[HTML]{FDCA7D}-0.4\% & \cellcolor[HTML]{EDE683}0.8\%  & \cellcolor[HTML]{EAE583}0.9\%  & \cellcolor[HTML]{E5E483}1.0\%  \\
\multicolumn{2}{c}{\multirow{-5}{*}{SP\_CTE4}}                  & 4                          & \cellcolor[HTML]{E9E583}0.9\%  & \cellcolor[HTML]{FDD47F}-0.1\% & \cellcolor[HTML]{C1DA81}1.8\%  & \cellcolor[HTML]{FDC57C}-0.5\% & \cellcolor[HTML]{FEE883}0.4\%  & \cellcolor[HTML]{FBEA84}0.5\%  & \cellcolor[HTML]{E2E383}1.1\%  & \cellcolor[HTML]{FCC27C}-0.6\% & \cellcolor[HTML]{FEE683}0.3\%  & \cellcolor[HTML]{CBDC81}1.6\%  \\ \hline
\end{tabular}
}
\caption{Step-wise and continuous $\alpha$-profiles compared with constant $\alpha$-profiles. Values represent the profit gap \% between the two, with respect to the constant profiles.}
\label{tab:cte_vs_cnt}
\end{table}

\begin{table}[t]
\centering
\resizebox{\textwidth}{!}{
\renewcommand\arraystretch{1.5}
\begin{tabular}{cccccccccccccccccc}
\hline
\multicolumn{18}{c}{Profit}                                                                                                                                                                                                                                                                                                                                                                                                                                                                                                                                                            \\ \hline
Instance & Myopic                            & SP\_CTE1                      & SP\_CTE2                       & SP\_CTE3                       & SP\_CTE4                       & SP\_CTE5                      & SP\_CTE6                       & SP\_STP1                       & SP\_STP2                       & SP\_STP3                       & SP\_STP4                       & SP\_STP5                       & SP\_STP6                       & SP\_PLY1                       & SP\_PLY2                       & SP\_PLY3                       & SP\_PLY4                       \\ \hline
0        & \cellcolor[HTML]{FA8E72}3076  & \cellcolor[HTML]{FDD47F}3255  & \cellcolor[HTML]{AAD380}3497   & \cellcolor[HTML]{B5D680}3474   & \cellcolor[HTML]{D1DE82}3413   & \cellcolor[HTML]{FEE883}3308  & \cellcolor[HTML]{F98871}3062   & \cellcolor[HTML]{C3DA81}3444   & \cellcolor[HTML]{AAD380}3498   & \cellcolor[HTML]{A8D27F}3501   & \cellcolor[HTML]{B1D580}3481   & \cellcolor[HTML]{A3D17F}3511   & \cellcolor[HTML]{C0D981}3449   & \cellcolor[HTML]{BBD881}3461   & \cellcolor[HTML]{A8D27F}3502   & \cellcolor[HTML]{95CD7E}3542   & \cellcolor[HTML]{BCD881}3458   \\
1        & \cellcolor[HTML]{FDD17F}3249  & \cellcolor[HTML]{FEE182}3290  & \cellcolor[HTML]{99CE7F}3534   & \cellcolor[HTML]{7FC67D}3590   & \cellcolor[HTML]{95CD7E}3542   & \cellcolor[HTML]{D3DF82}3409  & \cellcolor[HTML]{FBA777}3141   & \cellcolor[HTML]{89C97E}3568   & \cellcolor[HTML]{81C77D}3584   & \cellcolor[HTML]{B0D480}3485   & \cellcolor[HTML]{80C77D}3588   & \cellcolor[HTML]{93CC7E}3547   & \cellcolor[HTML]{A1D07F}3517   & \cellcolor[HTML]{A2D17F}3514   & \cellcolor[HTML]{87C97E}3571   & \cellcolor[HTML]{8CCA7E}3562   & \cellcolor[HTML]{9BCE7F}3530   \\
2        & \cellcolor[HTML]{FDD57F}3260  & \cellcolor[HTML]{F9EA84}3328  & \cellcolor[HTML]{8BCA7E}3563   & \cellcolor[HTML]{85C87D}3577   & \cellcolor[HTML]{8BCA7E}3564   & \cellcolor[HTML]{D7E082}3401  & \cellcolor[HTML]{FDC77D}3222   & \cellcolor[HTML]{7EC67D}3591   & \cellcolor[HTML]{73C37C}3614   & \cellcolor[HTML]{88C97E}3570   & \cellcolor[HTML]{89C97E}3568   & \cellcolor[HTML]{82C77D}3582   & \cellcolor[HTML]{6EC27C}3625   & \cellcolor[HTML]{96CD7E}3540   & \cellcolor[HTML]{63BE7B}3648   & \cellcolor[HTML]{87C97E}3572   & \cellcolor[HTML]{85C87D}3576   \\
3        & \cellcolor[HTML]{FCC47C}3216  & \cellcolor[HTML]{FEDD81}3280  & \cellcolor[HTML]{C3DA81}3443   & \cellcolor[HTML]{A8D27F}3502   & \cellcolor[HTML]{A2D07F}3515   & \cellcolor[HTML]{FEDC81}3278  & \cellcolor[HTML]{FCBC7B}3195   & \cellcolor[HTML]{98CE7F}3536   & \cellcolor[HTML]{ACD480}3492   & \cellcolor[HTML]{B7D680}3470   & \cellcolor[HTML]{90CB7E}3552   & \cellcolor[HTML]{8ECB7E}3556   & \cellcolor[HTML]{B7D680}3470   & \cellcolor[HTML]{A8D27F}3501   & \cellcolor[HTML]{94CD7E}3544   & \cellcolor[HTML]{93CC7E}3546   & \cellcolor[HTML]{91CC7E}3550   \\
4        & \cellcolor[HTML]{FA9373}3089  & \cellcolor[HTML]{FBAA77}3148  & \cellcolor[HTML]{F1E784}3346   & \cellcolor[HTML]{C9DC81}3430   & \cellcolor[HTML]{D7E082}3400   & \cellcolor[HTML]{FEE182}3290  & \cellcolor[HTML]{F8696B}2980   & \cellcolor[HTML]{C9DC81}3431   & \cellcolor[HTML]{DAE182}3395   & \cellcolor[HTML]{BBD881}3460   & \cellcolor[HTML]{E0E283}3382   & \cellcolor[HTML]{D2DE82}3412   & \cellcolor[HTML]{CFDD82}3418   & \cellcolor[HTML]{C7DB81}3436   & \cellcolor[HTML]{E1E383}3379   & \cellcolor[HTML]{D2DE82}3411   & \cellcolor[HTML]{BFD981}3453   \\ \hline
\multicolumn{18}{c}{Profit Gap \% with respect to the Myopic ILP}                                                                                                                                                                                                                                                                                                                                                                                                                                                                                                               \\ \hline
Instance & Myopic                            & SP\_CTE1                      & SP\_CTE2                       & SP\_CTE3                       & SP\_CTE4                       & SP\_CTE5                      & SP\_CTE6                       & SP\_STP1                       & SP\_STP2                       & SP\_STP3                       & SP\_STP4                       & SP\_STP5                       & SP\_STP6                       & SP\_PLY1                       & SP\_PLY2                       & SP\_PLY3                       & SP\_PLY4                       \\ \hline
0        & \cellcolor[HTML]{FA9A74}0.0\% & \cellcolor[HTML]{FFEB84}5.8\% & \cellcolor[HTML]{7CC67D}13.7\% & \cellcolor[HTML]{88C97E}12.9\% & \cellcolor[HTML]{AAD380}11.0\% & \cellcolor[HTML]{E3E383}7.5\% & \cellcolor[HTML]{FA9373}-0.5\% & \cellcolor[HTML]{99CE7F}12.0\% & \cellcolor[HTML]{7BC57D}13.7\% & \cellcolor[HTML]{7AC57D}13.8\% & \cellcolor[HTML]{85C87D}13.2\% & \cellcolor[HTML]{74C37C}14.1\% & \cellcolor[HTML]{96CD7E}12.1\% & \cellcolor[HTML]{8FCB7E}12.5\% & \cellcolor[HTML]{79C57D}13.8\% & \cellcolor[HTML]{63BE7B}15.1\% & \cellcolor[HTML]{91CC7E}12.4\% \\
1        & \cellcolor[HTML]{FA9A74}0.0\% & \cellcolor[HTML]{FBAB77}1.3\% & \cellcolor[HTML]{CEDD82}8.8\%  & \cellcolor[HTML]{B1D580}10.5\% & \cellcolor[HTML]{CADC81}9.0\%  & \cellcolor[HTML]{FEDE81}4.9\% & \cellcolor[HTML]{F86B6B}-3.3\% & \cellcolor[HTML]{BDD881}9.8\%  & \cellcolor[HTML]{B4D680}10.3\% & \cellcolor[HTML]{E7E483}7.3\%  & \cellcolor[HTML]{B2D580}10.4\% & \cellcolor[HTML]{C7DB81}9.2\%  & \cellcolor[HTML]{D7E082}8.2\%  & \cellcolor[HTML]{D8E082}8.2\%  & \cellcolor[HTML]{BBD881}9.9\%  & \cellcolor[HTML]{C0D981}9.6\%  & \cellcolor[HTML]{D0DE82}8.6\%  \\
2        & \cellcolor[HTML]{FA9A74}0.0\% & \cellcolor[HTML]{FCB77A}2.1\% & \cellcolor[HTML]{C5DB81}9.3\%  & \cellcolor[HTML]{BED981}9.7\%  & \cellcolor[HTML]{C5DB81}9.3\%  & \cellcolor[HTML]{FDD680}4.3\% & \cellcolor[HTML]{F98971}-1.2\% & \cellcolor[HTML]{B7D780}10.2\% & \cellcolor[HTML]{ABD380}10.9\% & \cellcolor[HTML]{C2DA81}9.5\%  & \cellcolor[HTML]{C3DA81}9.4\%  & \cellcolor[HTML]{BCD881}9.9\%  & \cellcolor[HTML]{A6D27F}11.2\% & \cellcolor[HTML]{D1DE82}8.6\%  & \cellcolor[HTML]{9ACE7F}11.9\% & \cellcolor[HTML]{C1D981}9.6\%  & \cellcolor[HTML]{BFD981}9.7\%  \\
3        & \cellcolor[HTML]{FA9A74}0.0\% & \cellcolor[HTML]{FCB579}2.0\% & \cellcolor[HTML]{EBE583}7.1\%  & \cellcolor[HTML]{CCDD82}8.9\%  & \cellcolor[HTML]{C5DB81}9.3\%  & \cellcolor[HTML]{FCB479}1.9\% & \cellcolor[HTML]{FA9172}-0.7\% & \cellcolor[HTML]{BAD881}10.0\% & \cellcolor[HTML]{D1DE82}8.6\%  & \cellcolor[HTML]{DDE182}7.9\%  & \cellcolor[HTML]{B2D580}10.4\% & \cellcolor[HTML]{B0D580}10.6\% & \cellcolor[HTML]{DDE182}7.9\%  & \cellcolor[HTML]{CDDD82}8.9\%  & \cellcolor[HTML]{B6D680}10.2\% & \cellcolor[HTML]{B5D680}10.3\% & \cellcolor[HTML]{B3D580}10.4\% \\
4        & \cellcolor[HTML]{FA9A74}0.0\% & \cellcolor[HTML]{FCB479}1.9\% & \cellcolor[HTML]{D6DF82}8.3\%  & \cellcolor[HTML]{A8D27F}11.0\% & \cellcolor[HTML]{B8D780}10.1\% & \cellcolor[HTML]{F4E884}6.5\% & \cellcolor[HTML]{F8696B}-3.5\% & \cellcolor[HTML]{A8D27F}11.1\% & \cellcolor[HTML]{BBD881}9.9\%  & \cellcolor[HTML]{98CE7F}12.0\% & \cellcolor[HTML]{C2DA81}9.5\%  & \cellcolor[HTML]{B2D580}10.5\% & \cellcolor[HTML]{AFD480}10.7\% & \cellcolor[HTML]{A5D17F}11.2\% & \cellcolor[HTML]{C4DA81}9.4\%  & \cellcolor[HTML]{B2D580}10.4\% & \cellcolor[HTML]{9CCF7F}11.8\% \\ \hline
\end{tabular}
}
\caption{Heat map of cumulative profit and profit gap (\%) w.r.t. Myopic ILP policy}
\label{tab:htmap}
\end{table}

\subsection{Profit acquired per interval}

\begin{figure}
  \centering
  \begin{subfigure}{0.6\textwidth}
    \includegraphics[width=\linewidth]{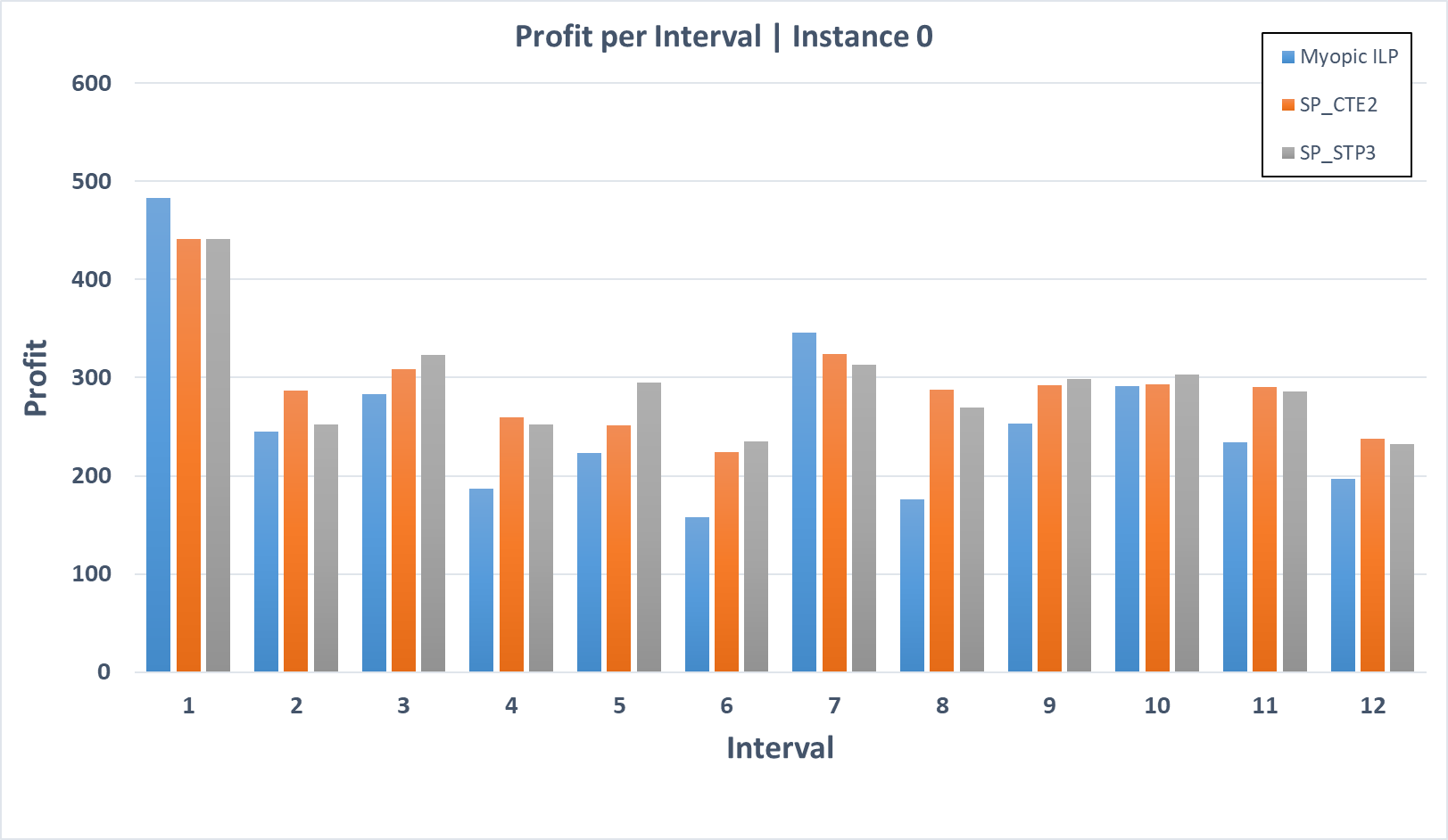}    
    \caption{}    
    \label{fig:subfig-a}
  \end{subfigure}
  \hfill
  \begin{subfigure}{0.6\textwidth}
    \includegraphics[width=\linewidth]{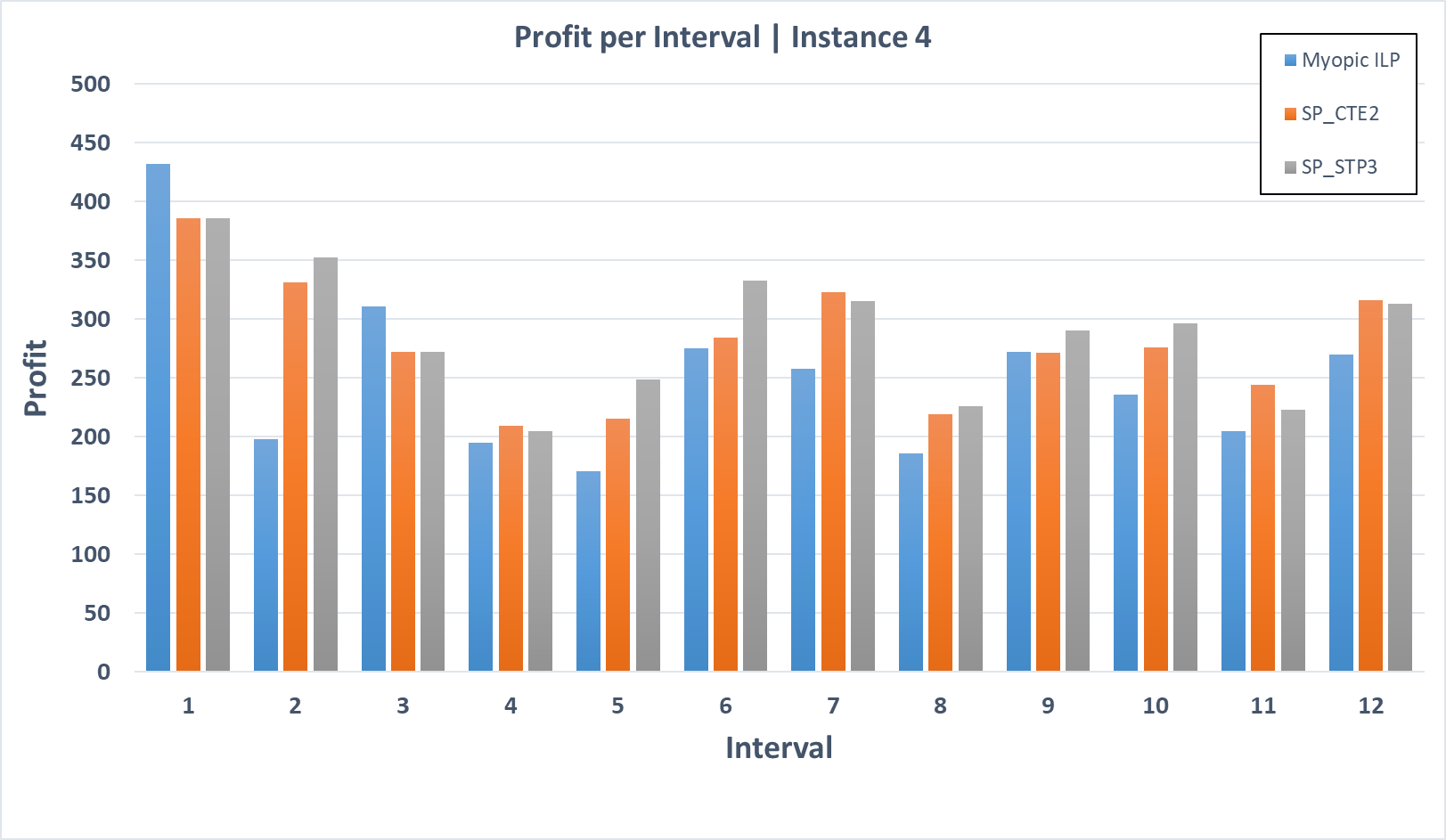}  
    \caption{}  
    \label{fig:subfig-b}
  \end{subfigure}
  \caption{Profit acquired per interval for Myopic ILP, SP\_CTE2 and SP\_STP3 in (a) Instance 0 and (b) Instance 4}
  \label{fig:combined}
\end{figure}

Figure~\ref{fig:combined} depicts the profit acquired per interval for test instances 0 and 4 by the Myopic ILP, SP\_CTE2 and SP\_STP3. In both instances the first interval witnessed a larger profit acquisition by the Myopic ILP. Thereafter the surrogate models closed the distance and maintained the performance gap. Presumably, in the first interval, link priority based slack creation limits the number of profitable requests being serviced. But once this reserve capacity is capitalized upon in subsequent intervals, the profit gained per interval is larger for the surrogate models. Occasionally the Myopic ILP policy performs slightly better (interval 7 of instance 0 and interval 3 of instance 4), but by and large the Surrogate ILP policies dominate and acquire a larger cumulative profit.
\\
\\
As denoted in Table \ref{tab:surr_legend}, SP\_CTE2 sets $\alpha_{i} = 1.5 \ \forall i \in I$ whereas 
SP\_STP3 transitions from 1.5 to 1.25 to 1 in steps of 4 intervals each. This means that for the first four intervals, both surrogate policies have identical $\alpha_{i}$ values ($= 1.5$). However the profit acquired per interval is different from the second interval onwards, owing to the dynamic and adaptive nature of the $k$NN component of the surrogate ILP policy.

\subsection{Significance of $\beta_{i,l}$}

We conclude the numerical experiments by demonstrating the role of using data-driven predictive prescriptions to determine efficient link priority vectors $\bm{\vec{\beta_{i}}} = \begin{bmatrix}
    \beta_{i,l}
\end{bmatrix}^{l \in A}$. The value of these estimations of link priority can be exhibited by comparison to an arbitrary policy which assumes equal importance (or priority) of all links in the network. Such a policy would assign a uniform link priority across the network. It would thereby delegate all responsibility for managing the differential slack in the network to the $\alpha$-profile. Yet it would suffer from the disadvantage of being unable to create reserve capacity differentiated across space (links in the network), and would only do this across time (varying $\alpha_{i}$ across intervals). To this end, we conducted experiments using the Surrogate ILP policy with a uniform standardized link priority vector $\bm{\vec{\beta_{i}}} = \begin{bmatrix}
    \frac{1}{|A|}
\end{bmatrix}^{l \in A}$ and report the outcome in Table~\ref{tab:beta_sig}. We observed that upon using a uniform link priority vector, the performance dropped sharply and in many cases even fell below the Myopic ILP policy benchmark. This shows that uniform slack creation across all network links fails to improve performance because the reserve capacities are introduced in a manner which cannot be taken advantage of. It is crucial to create reserve capacity differentiated across both space and time to outperform the benchmark of the myopic ILP. 

\begin{table}[t]
\centering
\resizebox{\textwidth}{!}{
\renewcommand\arraystretch{1.5}
\begin{tabular}{cccccccccccc}
\hline
\multirow{2}{*}{Instance} & \multirow{2}{*}{Total Arrived} & \multicolumn{2}{c}{Myopic ILP}                & \multirow{2}{*}{$\alpha$-profile} & \multirow{2}{*}{$\beta$ Application} & \multirow{2}{*}{Service Rate} & \multirow{2}{*}{Profit} & \multirow{2}{*}{Profit Gap} & \multirow{2}{*}{Profit Gap \%} & \multirow{2}{*}{Service Gap} & \multirow{2}{*}{Service Gap \%} \\ \cline{3-4}
                          &                                & \multicolumn{1}{l}{Service Rate} & Profit                &                            &                                      &                               &                         &                             &                                &                              &                                 \\ \hline
\multirow{6}{*}{0}        & \multirow{6}{*}{1203}          & \multirow{6}{*}{41.6\%}          & \multirow{6}{*}{3076} & \multirow{2}{*}{SP\_CTE2}  & kNN                                  & 40.6\%                        & 3497                    & 421                         & 13.7\%                         & -13                          & -1.0\%                          \\
                          &                                &                                  &                       &                            & uniform                              & 41.5\%                        & 3109                    & 33                          & 1.1\%                          & -2                           & -0.1\%                          \\
                          &                                &                                  &                       & \multirow{2}{*}{SP\_STP3}  & kNN                                  & 41.3\%                        & 3501                    & 425                         & 13.8\%                         & -4                           & -0.3\%                          \\
                          &                                &                                  &                       &                            & uniform                              & 42.5\%                        & 3157                    & 81                          & 2.6\%                          & 10                           & 0.9\%                           \\
                          &                                &                                  &                       & \multirow{2}{*}{SP\_PLY2}  & kNN                                  & 40.6\%                        & 3502                    & 426                         & 13.8\%                         & -13                          & -1.0\%                          \\
                          &                                &                                  &                       &                            & uniform                              & 41.3\%                        & 3091                    & 15                          & 0.5\%                          & -4                           & -0.3\%                          \\ \hline
\multirow{6}{*}{1}        & \multirow{6}{*}{1125}          & \multirow{6}{*}{45.5\%}          & \multirow{6}{*}{3249} & \multirow{2}{*}{SP\_CTE2}  & kNN                                  & 43.8\%                        & 3534                    & 285                         & 8.8\%                          & -19                          & -1.7\%                          \\
                          &                                &                                  &                       &                            & uniform                              & 46.4\%                        & 3244                    & -5                          & -0.2\%                         & 10                           & 0.9\%                           \\
                          &                                &                                  &                       & \multirow{2}{*}{SP\_STP3}  & kNN                                  & 43.6\%                        & 3485                    & 236                         & 7.3\%                          & -21                          & -1.9\%                          \\
                          &                                &                                  &                       &                            & uniform                              & 46.0\%                        & 3219                    & -30                         & -0.9\%                         & 6                            & 0.5\%                           \\
                          &                                &                                  &                       & \multirow{2}{*}{SP\_PLY2}  & kNN                                  & 44.3\%                        & 3571                    & 322                         & 9.9\%                          & -14                          & -1.2\%                          \\
                          &                                &                                  &                       &                            & uniform                              & 46.4\%                        & 3272                    & 23                          & 0.7\%                          & 10                           & 0.9\%                           \\ \hline
\multirow{6}{*}{2}        & \multirow{6}{*}{1260}          & \multirow{6}{*}{41.7\%}          & \multirow{6}{*}{3260} & \multirow{2}{*}{SP\_CTE2}  & kNN                                  & 39.7\%                        & 3563                    & 303                         & 9.3\%                          & -26                          & -2.0\%                          \\
                          &                                &                                  &                       &                            & uniform                              & 42.3\%                        & 3303                    & 43                          & 1.3\%                          & 7                            & 0.6\%                           \\
                          &                                &                                  &                       & \multirow{2}{*}{SP\_STP3}  & kNN                                  & 40.8\%                        & 3570                    & 310                         & 9.5\%                          & -12                          & -0.9\%                          \\
                          &                                &                                  &                       &                            & uniform                              & 41.7\%                        & 3313                    & 53                          & 1.6\%                          & -1                           & 0.0\%                           \\
                          &                                &                                  &                       & \multirow{2}{*}{SP\_PLY2}  & kNN                                  & 40.9\%                        & 3648                    & 388                         & 11.9\%                         & -11                          & -0.8\%                          \\
                          &                                &                                  &                       &                            & uniform                              & 40.7\%                        & 3221                    & -39                         & -1.2\%                         & -13                          & -1.0\%                          \\ \hline
\multirow{6}{*}{3}        & \multirow{6}{*}{1184}          & \multirow{6}{*}{44.5\%}          & \multirow{6}{*}{3216} & \multirow{2}{*}{SP\_CTE2}  & kNN                                  & 40.7\%                        & 3443                    & 227                         & 7.1\%                          & -45                          & -3.8\%                          \\
                          &                                &                                  &                       &                            & uniform                              & 42.6\%                        & 3091                    & -125                        & -3.9\%                         & -23                          & -1.9\%                          \\
                          &                                &                                  &                       & \multirow{2}{*}{SP\_STP3}  & kNN                                  & 42.1\%                        & 3470                    & 254                         & 7.9\%                          & -29                          & -2.4\%                          \\
                          &                                &                                  &                       &                            & uniform                              & 44.0\%                        & 3168                    & -48                         & -1.5\%                         & -6                           & -0.5\%                          \\
                          &                                &                                  &                       & \multirow{2}{*}{SP\_PLY2}  & kNN                                  & 42.2\%                        & 3544                    & 328                         & 10.2\%                         & -27                          & -2.3\%                          \\
                          &                                &                                  &                       &                            & uniform                              & 43.9\%                        & 3185                    & -31                         & -1.0\%                         & -7                           & -0.6\%                          \\ \hline
\multirow{6}{*}{4}        & \multirow{6}{*}{1147}          & \multirow{6}{*}{44.0\%}          & \multirow{6}{*}{3089} & \multirow{2}{*}{SP\_CTE2}  & kNN                                  & 41.0\%                        & 3346                    & 257                         & 8.3\%                          & -35                          & -3.0\%                          \\
                          &                                &                                  &                       &                            & uniform                              & 44.4\%                        & 3005                    & -84                         & -2.7\%                         & 4                            & 0.4\%                           \\
                          &                                &                                  &                       & \multirow{2}{*}{SP\_STP3}  & kNN                                  & 42.9\%                        & 3460                    & 371                         & 12.0\%                         & -13                          & -1.1\%                          \\
                          &                                &                                  &                       &                            & uniform                              & 44.0\%                        & 2975                    & -114                        & -3.7\%                         & 0                            & 0.0\%                           \\
                          &                                &                                  &                       & \multirow{2}{*}{SP\_PLY2}  & kNN                                  & 41.8\%                        & 3379                    & 290                         & 9.4\%                          & -26                          & -2.2\%                          \\
                          &                                &                                  &                       &                            & uniform                              & 43.8\%                        & 2976                    & -113                        & -3.7\%                         & -3                           & -0.2\%                          \\ \hline
\end{tabular}
}
\caption{Performance comparison with uniform $\beta$}
\label{tab:beta_sig}
\end{table}

\section{Conclusion and perspectives}
\label{con}

This study addresses the stochastic online version of the drone delivery service planning problem in urban airspace. The solution methods we explore are executed from the perspective of a network manager. Given an incoming Poisson process of drone delivery trip requests, we aim to find the optimal assignment of drones to space-time trajectories in an urban air traffic network within a pre-defined time horizon. We assume that drone deliveries must respect delivery time windows and that UAVs must fly at constant speed throughout the network to avoid airborne delays and congestion. Our main contribution is to expand the integer linear programming formulation of the deterministic DDSP to the stochastic variant with online demand. We also propose a novel solution methodology which leverages data-driven, machine learning-enabled predictions of network link priorities to construct a modified surrogate objective function. We outline a training scheme to train the machine learning component as well. The numerical results indicate that the Surrogate ILP policy outperforms the Myopic ILP policy which iteratively solves a for the objective of profit maximization. We demonstrate that even though the training procedure might be extensive, it is one-time investment returns rewards in the form of a sustained profit margin over the Myopic ILP policy across multiple instances. 
\\
\\
The escalating adoption of digital technology by the public and private sectors is producing real-time uncertain demand for logistic service providers globally. This is happening concurrently with the accelerating adoption of UAVs in transport fleets to service regions with varying degrees of urbanization. The utility of the solution methodology proposed in this study lies at the intersection of these two trends. By outperforming conventional online optimization methods our approach unlocks considerable benefits for the service provider and for the customers. Moreover, we successfully establish link priority and available link capacity as pivotal factors for enhancing network performance. This implies that the strategy and procedure of non-myopic anticipatory resource allocation we present in this paper is transferable to other classes of network-based logistic problems such as pickup and delivery logistics and ride-sharing. 
\\
\\
Multiple directions of further research can be identified. First, there is an opportunity to extend the data-driven optimization framework to learn efficient $\alpha$-profiles automatically. This also necessitates the design of a training scheme that explores the $\alpha$-profile space quickly and intelligently under time and computational budgets. Second, modifications can be imposed on the problem definition which make the DDSP more realistic, such as capacity constraints on each UAV, differential velocities and charging limitations. Extensions to the case of multiple deliveries per UAV also warrant further research. Although this stream of problem has been extensively studied through the lens of vehicle routing models, few studies have addressed the online demand case in a time-dependent routing context. The introduction and placement of charging platforms at selected nodes or a centralized drone hub somewhere in the network is also a compelling network-design possibility. 

\bibliography{reference}
\bibliographystyle{abbrvnat}

\end{document}